\documentclass[12pt,reqno]{amsart}
\usepackage{amsmath,amssymb,geometry,color,tikz-cd,extarrows}
\usepackage[backref,pagebackref,pdftex,hyperindex]{hyperref}
\usepackage[alphabetic]{amsrefs}
\usepackage{amssymb}
\usepackage{bbm}
\usepackage{enumitem}
\usepackage{upgreek}

\newtheorem{proposition}{Proposition}[section]
\newtheorem{theorem}[proposition]{Theorem}
\newtheorem{lemma}[proposition]{Lemma}
\newtheorem{corollary}[proposition]{Corollary}

\theoremstyle{definition}
\newtheorem{remark}[proposition]{Remark}
\newtheorem{definition}[proposition]{Definition}

\title{K-moduli with real coefficients}

\author{Yuchen Liu}
\address{Department of Mathematics, Northwestern University, Evanston, IL 60208, USA.}
\email{yuchenl@northwestern.edu}

\author{Chuyu Zhou}
\address{School of Mathematical Sciences, Xiamen University, Siming South Road 422, Xiamen, Fujian 361005, China.}
\email{chuyuzhou@xmu.edu.cn, chuyuzhou1@gmail.com}

\date{} 

\thanks{
        \newline
		\indent
		\emph{Competing interests}:  The authors have no conflict of interest to declare.
		}

\newcommand{\Fut}{{\rm{Fut}}}

\newcommand{\ord}{{\rm {ord}}}
\newcommand{\tc}{{\rm {tc}}}
\newcommand{\vol}{{\rm {vol}}}
\newcommand{\lct}{{\rm {lct}}}

\newcommand{\length}{{\rm {length}}}

\newcommand{\Spec}{{\rm {Spec}}}
\newcommand{\Proj}{{\rm{Proj}}}
\newcommand{\Exc}{{\rm {Exc}}}
\newcommand{\Val}{{\rm {Val}}}

\newcommand{\gr}{{\rm {gr}}}

\newcommand{\Supp}{{\rm {Supp}}}

\newcommand{\ac}{{\rm {ac}}}

\newcommand{\Hilb}{{\rm {Hilb}}}

\newcommand{\PGL}{{\rm {PGL}}}

\newcommand{\Aut}{{\rm {Aut}}}

\newcommand{\Ding}{{\rm {Ding}}}

\newcommand{\Kss}{{\rm {Kss}}}

\newcommand{\LF}{{\rm {LF}}}

\newcommand{\Kps}{{\rm {Kps}}}

\newcommand{\rank}{{\rm {rank}}}
\newcommand{\triv}{{\rm {triv}}}
\newcommand{\QM}{{\rm {QM}}}
\newcommand{\Int}{{\rm {Int}}}
\newcommand{\LCP}{{\rm {LCP}}}
\newcommand{\DivVal}{{\rm {DivVal}}}
\newcommand{\id}{{\rm {id}}}
\newcommand{\NA}{{\rm {NA}}}
\newcommand{\ST}{{\rm {ST}}}
\newcommand{\LEnv}{{\rm {LEnv}}}

\newcommand{\bA}{\mathbb{A}}

\newcommand{\bC}{\mathbb{C}}

\newcommand{\bG}{\mathbb{G}}

\newcommand{\bN}{\mathbb{N}}

\newcommand{\bP}{\mathbb{P}}
\newcommand{\bQ}{\mathbb{Q}}
\newcommand{\bR}{\mathbb{R}}

\newcommand{\bT}{\mathbb{T}}

\newcommand{\bZ}{\mathbb{Z}}

\newcommand{\mA}{\mathcal{A}}
\newcommand{\mB}{\mathcal{B}}

\newcommand{\mD}{\mathcal{D}}
\newcommand{\mE}{\mathcal{E}}
\newcommand{\mF}{\mathcal{F}}
\newcommand{\mG}{\mathcal{G}}

\newcommand{\mI}{\mathcal{I}}

\newcommand{\mK}{\mathcal{K}}
\newcommand{\mL}{\mathcal{L}}
\newcommand{\mM}{\mathcal{M}}

\newcommand{\mO}{\mathcal{O}}
\newcommand{\mP}{\mathcal{P}}

\newcommand{\mR}{\mathcal{R}}
\newcommand{\mS}{\mathcal{S}}

\newcommand{\mX}{\mathcal{X}}
\newcommand{\mY}{\mathcal{Y}}

\newcommand{\fB}{\mathbf{B}}

\newcommand{\fD}{\mathbf{D}}

\newcommand{\fL}{\mathbf{L}}

\newcommand{\fX}{\mathbf{X}}
\newcommand{\fY}{\mathbf{Y}}
\newcommand{\fZ}{\mathbf{Z}}

\newcommand{\fa}{\mathbf{a}}

\newcommand{\fr}{\mathbf{r}}

\newcommand{\tB}{\tilde{B}}

\newcommand{\tX}{\tilde{X}}

\newcommand{\ka}{\mathfrak{a}}
\newcommand{\kb}{\mathfrak{b}}

\begin{document}

\begin{abstract}
In this paper, we develop an algebraic K-stability theory (e.g. special test configuration theory and optimal destabilization theory) for log Fano $\bR$-pairs, and construct a proper K-moduli space to parametrize K-polystable log Fano $\bR$-pairs with some fixed invariants (e.g. dimension, volume, coefficients). All of these are well-known for log Fano $\bQ$-pairs, and the strategy in this paper is to reduce the problems (in many cases) to $\bQ$-coefficients case rather than rebuilding the whole program as in $\bQ$-coefficients case.

\end{abstract}

\maketitle

\setcounter{tocdepth}{1}

\tableofcontents

We work over the complex number field $\bC$ throughout the article.

\section{Introduction}

Based on a solid understanding of K-stability in algebraic sense (e.g. \cites{Oda13, LX14, Fuj19, Li17, LX20, LX18} etc.), a complete K-moduli theory is finally established for K-semistable log Fano $\bQ$-pairs (e.g. \cites{Jiang20, Xu20, BLX22, BX19, LWX21, BHLLX21, ABHLX20, XZ20b, CP21, LXZ22}). More precisely, there exists a finite type Artin stack parametrizing K-semistable log Fano $\bQ$-pairs with fixed invariants (e.g. dimension, volume, and Weil index), and the stack admits a projective and separated scheme as a good moduli space parametrizing K-polystable objects. See \cites{Xu21, Xu24} for a comprehensive understanding on this topic. 

Although it is expected that there is a parallel K-moduli theory for K-semistable log Fano $\bR$-pairs (for which we allow $\bR$-coefficients), it is seldom touched in this respect in the literature. One of the motivations of this article is to bridge this gap on the absence of the study for the K-stability of log Fano $\bR$-pairs. The natural idea is to check the whole algebraic K-stability theory and the whole K-moduli theory (in a bookkeeping sense) to see whether they apply for log Fano $\bR$-pairs. This is however a huge task with no doubt, though it is believed that many aspects should apply to the $\bR$-coefficients case without much difficulty. Our strategy in this article is to reduce the problems (in many cases) to $\bQ$-coefficients case and apply an approximation process, rather than duplicating the whole program.

To define the K-semistability of a log Fano $\bR$-pair, we directly take Fujita-Li's criterion as definition (see Definition \ref{def: Fujita-Li}). By an approximation process, we will see that the beta-invariant coincides with the generalized Futaki invariant for a special test configuration (see Proposition \ref{prop: fut=beta}). Under this definition, this article could be divided into two parts. The first part is Section \ref{sec: LX} and Section \ref{sec: optimal deg}, where we prove two fundamental results in algebraic K-stability theory corresponding to log Fano $\bR$-pairs, i.e. special test configuration theory and optimal destabilization.

\begin{theorem}\label{mainthm: LX}{\rm (Theorem \ref{thm: LX}, Proposition \ref{prop: weak LX})}
Let $(X, \Delta)$ be a log Fano $\bR$-pair. Then $(X, \Delta)$ is K-semistable if and only if $\Fut(\mX, \Delta_{\tc})\geq 0$ for any special test configuration $(\mX, \Delta_{\tc})\to \bA^1$ of $(X, \Delta)$. Moreover, if $(X, \Delta)$ is K-semistable and $(\mX, \Delta_\tc)\to \bA^1$ is a weakly special test configuration of $(X, \Delta)$ with $\Fut(\mX, \Delta_\tc)=0$, then the test configuration is special.
\end{theorem}

\begin{theorem}\label{mainthm: optimal deg}{\rm (Theorem \ref{thm: special minimizer}, Theorem \ref{thm: optimal deg})}
Let $(X, \Delta)$ be a log Fano $\bR$-pair of dimension $d$ with $\delta(X, \Delta)< \frac{d+1}{d}$. Then there exists a divisorial valuation computing $\delta(X, \Delta)$, and any such divisorial valuation is special (i.e. induced by a special test configuration). Moreover, if $\delta(X, \Delta)\leq 1$, then any divisorial valuation computing $\delta(X, \Delta)$ is induced by a special test configuration whose central fiber has the same $\delta$-invariant as $(X, \Delta)$.
\end{theorem}

For log Fano $\bQ$-pairs, Theorem \ref{mainthm: LX} is known by \cite{LX14} and Theorem \ref{mainthm: optimal deg} is known by \cite{LXZ22}. However, for Theorem \ref{mainthm: LX}, it is weaker than the results in \cite{LX14} because we only work out the last step (of \cite{LX14}) of the process of decreasing the generalized Futaki invariants (e.g. Proposition \ref{prop: weak LX}); for Theorem \ref{mainthm: optimal deg}, the different point from \cite[Theorem 1.2]{LXZ22} is that the $\delta$-invariant is not necessarily rational here. 

The second part of this article consists of Sections \ref{sec: stack}, \ref{sec: Theta}, \ref{sec: S}, and  \ref{sec: proper}, where we develop a K-moduli theory for K-semistable log Fano $\bR$-pairs.
To state the result, we first define a set of log Fano pairs, denoted by $\mK:=\mK(d, v, I)$, where $d$ is a positive integer, $v$ is a positive real number, and $I$ is a finite set of positive real numbers. We say a log Fano pair $(X, \Delta)$ is contained in $\mK$ if and only if the following conditions are satisfied:
\begin{enumerate}
\item $(X, \Delta)$ is a K-semistable log Fano pair of dimension $d$;
\item the coefficients of $\Delta$ are contained in $I$;
\item $\vol(-K_X-\Delta)\geq v$.
\end{enumerate}
Fixing a real number $v_0\geq v$ and an irrational vector $\fa:=(a_1,...,a_k)$, where $a_j\in I$ for each $1\leq j\leq k$, we consider the following subset of $\mK$:
$$\mK_{d, v_0, \fa}:=\{(X, \sum_{j=1}^k a_jD_j)\in \mK\ |\ \textit{$D_j$'s are Weil divisors and $\vol(-K_X-\sum_{j=1}^k a_jD_j)=v_0$}\}. $$

\begin{theorem}\label{mainthm: moduli}{\rm (Theorem \ref{thm: artin stack}, Theorem \ref{thm: gms})}
There exists a finite type Artin stack $\mM^\Kss_{d, v_0, \fa}$ parametrizing log Fano $\bR$-pairs in $\mK_{d, v_0, \fa}$, and the stack admits a proper and separated good moduli space $M^\Kps_{d, v_0, \fa}$ parametrizing K-polystable objects.
\end{theorem}

Theorem \ref{mainthm: moduli} is established by verifying several main ingredients as in the case of log Fano $\bQ$-pairs, e.g. boundedness (Theorem \ref{thm: bdd}), Openness (Theorem \ref{thm: openness}), $\Theta$-reductivity (Theorem \ref{thm: Theta}), $S$-completeness (Theorem \ref{thm: S}), Properness (Theorem \ref{thm: proper R}).

The second motivation for this article is to present a more complete wall-crossing picture in the non-proportional wall crossing theory. In \cite{LZ24}, we formulate the non-proportional wall crossing problem for K-stability and work out a semi-algebraic chamber decomposition to control the variation of K-semi(poly)stability (e.g. \cite[Section 1.1]{LZ24}).  Armed with Theorem \ref{mainthm: moduli}, we are able to present a  wall-crossing diagram which takes into account of K-moduli with real coefficients. Please refer to \cite[Section 9]{LZ24}.

Finally, we expect that the K-moduli space in the $\bR$-coefficient case is projective with an ample CM $\bR$-line bundle, which would be a natural generalization of \cite{CP21, XZ20b}. This is not discussed in this paper.

\noindent
\subsection*{Acknowledgement}
YL is partially supported by NSF CAREER Grant DMS-2237139 and an AT\&T Research Fellowship from Northwestern University.
CZ is supported by an NSFC grant (No. 12501058), a grant from Xiamen University (No. X2450214), and Samsung Science and Technology Foundation under Project Number SSTF-BA2302-03. We thank Harold Blum for comments.

\section{Preliminaries}

In this section, we will review some preliminaries that will be applied in this article.   

We say $(X, \Delta)$ is a \textit{couple} if $X$ is a normal variety and $\Delta$ is an effective $\bR$-divisor on $X$. 
We say a couple $(X, \Delta)$ is a \textit{log pair} if $K_X+\Delta$ is $\bR$-Cartier. 
For a projective log pair $(X, \Delta)$, we say it is \textit{log Fano} if $(X, \Delta)$ admits klt singularities and $-K_X-\Delta$ is ample; in the case $\Delta=0$, we just say \textit{Fano} instead of log Fano.  For a projective normal variety $X$, we say it is Fano type if there exists an effective $\bR$-divisor $\Delta$ on $X$ such that $(X, \Delta)$ is log Fano.

Sometimes we want to emphasize the coefficients of a log Fano pair, so we make the following convention. For a log Fano pair $(X, \Delta)$, if $\Delta$ is a $\bQ$-divisor, we call it a \textit{log Fano $\bQ$-pair}; if $\Delta$ admits irrational coefficients, we call it a \textit{log Fano $\bR$-pair}. We highlight here that log Fano $\bR$-pairs always admit irrational coefficients by the convention. 

A \textit{log smooth model} $(Y, E)$ over a couple $(X, \Delta)$ consists of a log resolution $\pi: Y\to (X, \Delta)$ and a reduced divisor $E$ on $Y$ such that $E+\Exc(\pi)+\pi_*^{-1}\Delta$ has simple normal crossing (SNC) support.

For various types of singularities associated to a log pair,  e.g. \textit{Kawamata log terminal (klt), log canonical (lc)}, etc., we refer to \cite{KM98, Kollar13}.

\subsection{Valuations}

Let $X$ be a normal variety. A valuation $v$ on $X$ is an $\bR$-valued function $v: K(X)^*\to \bR$ (where $K(X)$ is the function field and $K(X)^*=K(X)\setminus 0$) satisfying: (1) $v(fg)=v(f)+v(g)$; (2) $v(f+g)\geq \min\{v(f), v(g)\}$; (3) $v|_{\bC^*}=0$. We set $v(0)=+\infty$. The \textit{trivial valuation $v_{\triv}$} is defined by $v_{\triv}(f)=0$ for any $f\in K(X)^*$. A \textit{divisorial valuation} is of the form $c\cdot\ord_F$, where $c\in \bR_{>0}$ and $F$ is a prime divisor over $X$.
We denote $\Val_X$ to be the set of all valuations on $X$, and write $\Val_X^*:=\Val_X\setminus v_\triv$.

\begin{definition}
Given a normal variety $X$. An \textit{ideal sequence} on $X$ is a sequence of ideals contained in $\mO_X$, denoted by $\ka_\bullet:=\{\ka_m\}_{m\in \bN}$, satisfying $\ka_0=\mO_X$ and $\ka_{m_1}\ka_{m_2}\subset \ka_{m_1+m_2}$ for any $m_1,m_2\in \bN$.  Given a valuation $v\in \Val_X$, for any $p\in \bR_{\geq 0}$ we set
$$\ka_p(v):=\{f\in \mO_X\ |\ \textit{$v(f)\geq p$}\}. $$
The \textit{ideal sequence} of $v$ is denoted by $\ka_\bullet(v):=\{\ka_m(v)\}_{m\in \bN}$.
\end{definition}

\begin{definition}\label{def: qm}
Given a normal variety $X$ and a birational morphism $Y\to X$ where $Y$ is normal. Let $\eta$ be a scheme-theoretic point such that $\mO_{Y, \eta}$ is regular and let $(y_1,...,y_r)$ be a regular system of parameters. For any $\alpha:=(\alpha_1,...,\alpha_r)\in \bR^r_{\geq 0}$, there is a natural way to define a valuation $v_\alpha\in \Val_X$. For any $f\in \mO_{Y, \eta}$, we could write $f=\sum_{\beta\in \bZ^r_{\geq 0}}c_\beta\cdot y^{\beta}$ in $\widehat{\mO_{Y, \eta}}\cong \kappa(\eta)[\![y_1,...,y_r]\!]$, where $c_\beta\in \kappa(\eta)$ and $y^\beta=y_1^{\beta_1}...y_r^{\beta_r}$ with $\beta=(\beta_1,...,\beta_r)\in \bZ^r_{\geq 0}$. We set
$$v_\alpha(f):=\min\{\langle\alpha, \beta\rangle\ |\ \textit{$c_\beta\ne 0$}\}. $$
It is not hard to check that $v_\alpha$ is a valuation and such a valuation is called a \textit{quasi-monomial valuation} on $X$.

If, in addition, $(Y, E:=\sum_{i=1}^lE_i)\to X$ is a log smooth model 
where $(y_i=0)=E_i$ for each $1\leq i\leq r$ as an irreducible component of $E$, then we denote the set $\{v_\alpha\ |\ \textit{$\alpha\in \bR^r_{\geq 0}$}\}$ by $\QM_\eta(Y, E)$. Write $\QM(Y, E):=\cup_\eta \QM_\eta(Y, E)$, where $\eta$ runs over all generic points of strata of $E$.
\end{definition}

\begin{definition}
For a log pair $(X, \Delta)$, we define the \textit{log discrepancy} for a non-trivial valuation on $X$. For a divisorial valuation $v=c\cdot \ord_F$, we define
$$A_{X, \Delta}(v):=c\cdot A_{X, \Delta}(F):=c\cdot(\ord_{F}\left(K_Y-f^*(K_X+\Delta)\right)+1), $$  
where $f: Y\to X$ is a normal birational model such that $F$ is a prime divisor on $Y$.
For a quasi-monomial valuation $v_\alpha$ as in Definition \ref{def: qm}, we define
$$A_{X, \Delta}(v_\alpha):=\sum_{i=1}^r\alpha_i\cdot A_{X, \Delta}(E_i). $$
For the log discrepancy of a more general valuation, please refer to \cite{JM12}. We denote $\Val_X^\circ:=\{v\in \Val_X^*\ |\ \textit{$A_{X, \Delta}(v)<+\infty$}\}$. For any valuation $v\in \Val_X^\circ$ and any $c\in \bR_{>0}$, we have $A_{X, \Delta}(c\cdot v)=c\cdot A_{X, \Delta}(v)$.
\end{definition}

\subsection{$S$-invariants}\label{subsec: S-inv}

In this subsection, we fix $(X, \Delta)$ to be a log Fano $\bR$-pair of dimension $d$ (since all the concepts in this subsection are standard for log Fano $\bQ$-pairs, we only consider log Fano $\bR$-pairs). For a large positive integer $m\in \bZ_{>0}$, we say $D_m$ is an \textit{$m$-basis type divisor for $(X, \Delta)$} if it is of the following form
$$D_m=\frac{\sum_{j=1}^{N_m}{\rm div}(s_j)}{mN_m}\sim_\bR -(K_X+\Delta),  $$
where $N_m=\dim H^0(X,- m(K_X+\Delta))$ and $\{s_j\}_{j=1}^{N_m}$ is a basis of the vector space $H^0(X, - m(K_X+\Delta))$. Note the following convention
$$ H^0(X, - m(K_X+\Delta)):=H^0(X, \lfloor- m(K_X+\Delta)\rfloor)+\lceil m\Delta\rceil -m\Delta,$$
and we understand a section of $H^0(X, - m(K_X+\Delta))$ as a section of  $H^0(X, \lfloor- m(K_X+\Delta)\rfloor)$ plus the fraction part $\lceil m\Delta\rceil -m\Delta$.
Then we define the following invariant for a valuation $v\in \Val_X^\circ$:
$$S_m(v):=\sup_{D_m} v (D_m), $$ 
where $D_m$ runs through all $m$-basis type divisors. By the same arguments in \cite{AZ22, BJ20}, we know that $\lim_{m\to \infty} S_m(v)$ exists and denote it by $S_{X, \Delta}(v)$. Moreover, when $v=c\cdot \ord_F$ is a divisorial valuation, the limit of $S_m(v)=c\cdot S_m(F)$ can be reformulated as follows:
$$S_{X, \Delta}(v)=c\cdot S_{X, \Delta}(F)= \frac{c}{\vol(-K_X-\Delta)}\int_0^{\infty}\vol(-f^*(K_X+\Delta)-tF){\rm d} t,$$
where $f: Y\to X$ is a normal proper birational model such that $F$ is a prime divisor on $Y$.
For each large positive integer $m$, we define:
$$\delta_m(X, \Delta):=\inf_{E}\frac{A_{X, \Delta}(E)}{S_m(E)}, $$
where $E$ runs through all prime divisors over $X$. By the same arguments in \cite{AZ22, BJ20}, we know that $\lim_{m\to \infty}\delta_m(X, \Delta)$ exists and  
$$\delta(X, \Delta):=\lim_{m\to \infty}\delta_m(X, \Delta)=\inf_{v\in \Val_X^\circ}\frac{A_{X, \Delta}(v)}{S_{X, \Delta}(v)}=\inf_{E} \frac{A_{X, \Delta}(E)}{S_{X, \Delta}(E)}, $$
where $E$ runs through all prime divisors over $X$. If there exists a valuation  $v\in \Val_X^\circ$ (resp. a prime divisor $E$ over $X$) satisfying $\delta(X, \Delta)=\frac{A_{X, \Delta}(v)}{S_{X, \Delta}(v)}$ (resp. $\delta(X, \Delta)=\frac{A_{X, \Delta}(E)}{S_{X, \Delta}(E)}$), we say $v$ (resp. $E$) computes $\delta(X, \Delta)$.

Suppose $G$ and $D$ are effective $\bR$-divisors on $X$, we define 
$$\ord_{G} (D):=\sup\{\lambda\in \bR_{\geq 0}\ |\ \textit{$D\geq \lambda\cdot G$}\}.$$ 
Note that $\ord_G$ is not necessarily a valuation on $X$. Define $S_m(G):=\sup_{D_m} \ord_G(D_m)$, where $D_m$ runs through all $m$-basis type divisors. By \cite{AZ22}, we have
$$S_{X, \Delta}(G):=\lim_{m\to \infty} S_m(G)= \frac{1}{\vol(-K_X-\Delta)}\int_0^{\infty}\vol(-K_X-\Delta-tG){\rm d} t.$$
In particular, if $G\in |-K_X-\Delta|_\bR$,  a simple calculation implies $S_{X, \Delta}(G)=\frac{1}{d+1}$.

\begin{proposition}\label{prop: compatible}
Let $(X, \Delta)$ be a log Fano $\bR$-pair and $v\in \Val_X^\circ$ a valuation.  Suppose $G$ is an effective $\bR$-divisor on $X$. Then for each positive integer $m$ such that $H^0(X,-m(K_X+\Delta))\neq 0$, there exists an $m$-basis type divisor $D_m$ such that
$$S_m(v)=v(D_m)\quad \text{and}\quad S_m(G)=\ord_G(D_m). $$
\end{proposition}

\begin{proof}
Implied by \cite[Lemma 3.1]{AZ22}.
\end{proof}

\subsection{K-stability}

Inspired by \cite{Fuj19} and  \cite{Li17}, we define the K-semistability of a log Fano pair as follows.
\begin{definition} (Fujita--Li)\label{def: Fujita-Li}
For a given log Fano pair $(X, \Delta)$, we say it is \textit{K-semistable} if 
$\beta_{X, \Delta}(E):=A_{X, \Delta}(E)-S_{X, \Delta}(E)\geq 0$ for any prime divisor $E$ over $X$.
\end{definition}

The following result is clear given Definition \ref{def: Fujita-Li}.

\begin{theorem}{\rm (\cite{Fuj19, Li17, FO18, BJ20})}
Let $(X, \Delta)$ be a log Fano pair. Then $(X, \Delta)$ is K-semistable if and only if $\delta(X, \Delta)\geq 1$.
\end{theorem}

We also recall the concept of K-polystability. First we need the following concept of test configurations.

\begin{definition}{\rm (Test configuration)}\label{def: tc}
Let $(X,\Delta)$ be a projective log pair of dimension $d$ and $L$ an ample $\bR$-line bundle on $X$. A \textit{test configuration} (of $(X, \Delta; L)$) $\pi: (\mX,\Delta_\tc;\mL)\to \bA^1$ is a degenerating family over $\bA^1$ consisting of the following data:
\begin{enumerate}
\item $\pi: \mX\to \bA^1$ is a projective flat morphism of \textit{normal} varieties, $\Delta_\tc$ is an effective $\bR$-divisor on $\mX$, and $\mL$ is a relatively ample $\bR$-line bundle on $\mX$;
\item the family $\pi$ admits a $\bG_m$-action which lifts the standard $\bG_m$-action on $\bA^1$ such that $(\mX,\Delta_\tc; \mL)\times_{\bA^1}(\bA^1\setminus \{0\})$ is $\bG_m$-equivariantly isomorphic to $(X, \Delta; L)\times (\bA^1\setminus \{0\})$ where the $\bG_m$-action is trivial on the first component and standard on the second component.
\end{enumerate}
We say the test configuration $(\mX, \Delta_\tc;\mL)$ is \textit{product type} if it is induced by a one parameter subgroup of $\Aut(X, \Delta; L)$; we say the test configuration is \textit{trivial} if it is $\bG_m$-equivariantly isomorphic to $(X, \Delta; L)\times \bA^1$. 

We denote $(\bar{\mX}, \bar{\Delta}_\tc;\bar{\mL})\to \bP^1$ to be the natural compactification of the original test configuration, which is obtained by glueing $(\mX, \Delta_\tc;\mL)$ and $(X,\Delta;L)\times (\bP^1\setminus 0)$ along their common open subset $(X, \Delta;L)\times \bA^1\setminus\{0\}$.
\textit{The generalized Futaki invariant} associated to $(\mX, \Delta_\tc; \mL)$ is defined as follows:
$$\Fut(\mX,\Delta_\tc;\mL):=\frac{((K_{\bar{\mX}/\bP^1}+\bar{\Delta}_\tc)\cdot\bar{\mL}^d)}{(L^d)} -\frac{d}{d+1}\frac{\left((K_X+\Delta)L^{d-1}\right)\bar{\mL}^{d+1}}{(L^d)^2}.$$
\end{definition}

\begin{definition}{\rm (Special \& Weakly special test configurations)}\label{def: special tc}
Suppose $(X,\Delta)$ is a log Fano pair of dimension $d$ and $L=-K_X-\Delta$. Let  $(\mX,\Delta_\tc; \mL)$ be a test configuration of $(X, \Delta;L)$ such that $\mL=-K_{\mX/\bA^1}-\Delta_\tc$. We call it a \textit{weakly special (resp. special) test configuration} if  $(\mX, \mX_0+\Delta_{\tc})$ admits \textit{lc} (resp. \textit{plt}) singularities. Suppose $(\mX, \Delta_\tc;\mL)\to \bA^1$ is a non-trivial \textit{special} test configuration, then ${\ord_{\mX_0}}|_{K(X)}=c\cdot \ord_E$ for some $c\in \bQ_{>0}$, where $K(X)$ is the function field of $X$ and $E$ is a prime divisor over $X$ (e.g. \cite{BHJ17}), and we call $E$ (resp. $a\cdot \ord_E$ for any $a>0$) a \textit{special divisor (resp. special divisorial valuation)} for $(X, \Delta)$.
The generalized Futaki invariant of a weakly special test configuration is simplified as follows:
$$\Fut(\mX,\Delta_\tc;\mL)=-\frac{1}{d+1}\cdot\frac{(-K_{\bar{\mX}/\bP^1}-\bar{\Delta}_\tc)^{d+1}}{(-K_X-\Delta)^d}.$$
We will leave out the polarizations $L$ and $\mL$ when there is no confusion.
\end{definition}

For log Fano pairs, we do not consider more general test configurations (e.g. worse than weakly special ones) in this paper. The following definition is inspired by \cite{LWX21}.

\begin{definition}\label{def: kps}
Let $(X, \Delta)$ be a K-semistable log Fano pair. We say it is  \textit{K-polystable} if any special test configuration of $(X, \Delta; -K_X-\Delta)$ with a K-semistable central fiber is product type.
\end{definition}

These definitions (i.e. Definitions \ref{def: Fujita-Li}, \ref{def: tc}, \ref{def: special tc}, \ref{def: kps}) are standard in literature for log Fano $\bQ$-pairs (e.g. \cite{Xu21}). To study log Fano $\bR$-pairs, our philosophy is to approximate log Fano $\bR$-pairs by log Fano $\bQ$-pairs, and reduce the problem to the $\bQ$-coefficients case. The following continuity property is crucial to achieve the approximation process.

\begin{proposition}\label{prop: continuity}
Let $(X, \Delta)$ be a log Fano pair. Let $\{(X, \Delta_i)\}_{i\in \bN}$ be a sequence of log Fano $\bQ$-pairs such that each $\Delta_i$ supports in $\Supp (\Delta)$ and $||\Delta_i-\Delta||$ tends to zero as $i\to \infty$. Then $\lim_{i\to \infty} \delta(X, \Delta_i)=\delta(X, \Delta)$.
\end{proposition}

\begin{proof}
This argument is essentially proved in \cite[Theorem 1.7]{Zhang21}, where the setting is slightly different. In our setting, it suffices to show that for any $0<\epsilon< 1$, there exists a  positive integer $N$ satisfying the following $(\star)$
$$\left|\inf_{E}\frac{A_{X, \Delta_i}(E)}{S_{X, \Delta_i}(E)}-\inf_E\frac{A_{X, \Delta}(E)}{S_{X, \Delta}(E)}\right|<\epsilon $$
for any $i>N$, where $E$ runs over all prime divisors over $X$. 

Up to a small $\bQ$-factorization, we may assume $X$ is $\bQ$-factorial. For any fixed $0<\epsilon< 1$, there exists a positive integer $N$ such that $\epsilon (-K_X-\Delta)-|\Delta-\Delta_i|$
is $\bR$-linear equivalent to an effective $\bR$-divisor for $i>N$.
Thus we have the following $(\clubsuit)$
$$|A_{X, \Delta}(E)-A_{X, \Delta_i}(E)|=|\ord_E(\Delta_i-\Delta)|\leq \epsilon\cdot T_{X,\Delta}(E)$$
for any prime divisor $E$ over $X$ and any $i>N$, where $T_{X, \Delta}(E)$ is the pseudo-effective threshold of $E$ with respect to $-K_X-\Delta$.
On the other hand, there exists a positive integer $N'$ such that
$$(1-\epsilon)(-K_X-\Delta)\leq -K_X-\Delta_i\leq (1+\epsilon)(-K_X-\Delta) $$
for any $i>N'$. Thus we easily deduce the following $(\clubsuit\clubsuit)$
$$\frac{(1-\epsilon)^{d+1}}{(1+\epsilon)^d}\cdot S_{X, \Delta}(E)\leq S_{X, \Delta_i}(E)\leq \frac{(1+\epsilon)^{d+1}}{(1-\epsilon)^d}\cdot S_{X, \Delta}(E)$$
for any prime divisor $E$ over $X$ and any $i>N'$. Combining $(\clubsuit)$ and $(\clubsuit\clubsuit)$, we have
$$\inf_{E}\frac{A_{X, \Delta}(E)-\epsilon\cdot T_{X,\Delta}(E)}{\frac{(1+\epsilon)^{d+1}}{(1-\epsilon)^d}\cdot S_{X, \Delta}(E)} \leq \inf_{E}\frac{A_{X, \Delta_i}(E)}{S_{X, \Delta_i}(E)}\leq \inf_{E}\frac{A_{X, \Delta}(E)+\epsilon\cdot T_{X,\Delta}(E)}{\frac{(1-\epsilon)^{d+1}}{(1+\epsilon)^d}\cdot S_{X, \Delta}(E)} $$
for any $i>N+N'$.
Recalling the inequality $\frac{T_{X, \Delta}(E)}{S_{X, \Delta}(E)}\leq d+1$ (e.g. \cite{BJ20}),
we have
$$\inf_{E}\frac{A_{X, \Delta_i}(E)}{S_{X, \Delta_i}(E)}\leq \frac{(1+\epsilon)^{d}}{(1-\epsilon)^{d+1}}\left(\inf_{E}\frac{A_{X, \Delta}(E)}{S_{X, \Delta}(E)} + \epsilon\cdot (d+1)\right)$$
and 
$$\frac{(1-\epsilon)^{d}}{(1+\epsilon)^{d+1}}\left(\inf_{E}\frac{A_{X, \Delta}(E)}{S_{X, \Delta}(E)} - \epsilon\cdot (d+1)\right)\leq \inf_{E}\frac{A_{X, \Delta_i}(E)}{S_{X, \Delta_i}(E)} $$
for any $i>N+N'$. The above two inequalities imply $(\star)$. The proof is complete.
\end{proof}

\subsection{Approximation tools}

In this subsection, we collect some approximation tools which help to study log Fano $\bR$-pairs via approximation by log Fano $\bQ$-pairs. 

\begin{proposition}\label{prop: QR}
Let $(X, \Delta:=\sum_{j=1}^k r_j D_j)$ be a log Fano pair, where $D_j$'s are effective Weil divisors and $(r_1,...,r_k)\in \bR^k_{\geq 0}$ is  irrational. Suppose $(\mX, \Delta_\tc:=\sum_{j=1}^kr_j \mD_j)\to \bA^1$ is a special test configuration of $(X, \Delta)$. Then there exists a rational polytope $P\subset \bR^k_{\geq 0}$ containing $(r_1,...,r_k)$ as an interior point such that we have the following conclusions:
\begin{enumerate}
\item $(X, \sum_{j=1}^k x_j D_j)$ is a log Fano pair for any $(x_1,...,x_k)\in P$;
\item $(\mX, \sum_{j=1}^kx_j\mD_j)\to \bA^1$ is a special test configuration of $(X, \sum_{j=1}^kx_jD_j)$ for any $(x_1,...,x_k)\in P$.
\end{enumerate}
\end{proposition}

\begin{proof}
Implied by \cite[Corollary 7.4]{LZ24}.
\end{proof}

\begin{proposition}\label{prop: general QR}
Given a family of log Fano pairs 
$(Y, B:=\sum_{j=1}^k r_j B_j)\to T$ over a smooth base $T$, where $Y\to T$ is flat, $B_j$'s are effective Weil divisors, and $(r_1,...,r_k)\in \bR^k_{\geq 0}$ is irrational. Suppose $-K_{Y/T}-B$ is $\bR$-Cartier. Then there exists a rational polytope $P\subset \bR^k_{\geq 0}$ containing $(r_1,...,r_k)$ as an interior point such that we have the following conclusions:
\begin{enumerate}
\item $(Y,\sum_{j=1}^k x_j B_j)\to T$ is a flat family of log Fano pairs for any $(x_1,...,x_k)\in P$;
\item $-K_{Y/T}-\sum_{j=1}^kx_jB_j$ is $\bR$-Cartier for any $(x_1,...,x_k)\in P$.
\end{enumerate}
\end{proposition}

\begin{proof}
Implied by \cite[Proposition 7.3]{LZ24}.
\end{proof}

The following result is an application of Proposition \ref{prop: QR}.

\begin{proposition}\label{prop: fut=beta}
Let $(X, \Delta:=\sum_{j=1}^k r_j D_j)$ be a log Fano pair, where $D_j$'s are effective Weil divisors. Suppose $(\mX, \Delta_\tc:=\sum_{j=1}^kr_j \mD_j)\to \bA^1$ is a non-trivial special test configuration of $(X, \Delta)$ and write ${\ord_{\mX_0}}|_{K(X)}=c\cdot \ord_E$, where $K(X)$ is the function field of $X$ and $c\in \bQ_{>0}$. Then we have the following equality:
$$\Fut(\mX, \Delta_\tc)=c\cdot (A_{X, \Delta}(E)-S_{X, \Delta}(E)).$$ 
\end{proposition}

\begin{proof}
If $(r_1,...,r_k)$ is a rational vector, then the formula is well-known by \cite{Fuj19}. We may assume $(r_1,...,r_k)$ is irrational.
By Proposition \ref{prop: QR}, there exists a rational polytope $P$ containing $(r_1,...,r_k)$ as an interior point such that 
\begin{enumerate}
\item $(X, \sum_{j=1}^k x_j D_j)$ is a log Fano pair for any $(x_1,...,x_k)\in P$;
\item $(\mX, \sum_{j=1}^kx_j\mD_j)\to \bA^1$ is a special test configuration of $(X, \sum_{j=1}^kx_jD_j)$ for any $(x_1,...,x_k)\in P$.
\end{enumerate}
Choose a sequence of rational vectors $\{(r_{i1},...,r_{ik})\}_{i=1}^\infty\subset P$ tending to $(r_1,...,r_k)$.
By \cite{Fuj19}, we have the following formula for each $i$:
$$\Fut(\mX, \sum_{j=1}^kr_{ij}\mD_j; -K_{\mX}-\sum_{j=1}^kr_{ij}\mD_j) =c\cdot \left(A_{X, \sum_{j=1}^kr_{ij}D_j}(E)-S_{X, \sum_{j=1}^kr_{ij}D_j}(E)\right).$$
By the continuity of both sides, we have the desired formula
$$\Fut(\mX, \Delta_{\tc}; -K_{\mX}-\Delta_{\tc}) =c\cdot \left(A_{X, \Delta}(E)-S_{X, \Delta}(E)\right).$$
The proof is complete.
\end{proof}

\begin{corollary}\label{cor: kss deg}
Given a K-semistable log Fano pair $(X, \Delta)$. Suppose $(\mX, \Delta_\tc)\to \bA^1$ is a special test configuration such that the central fiber $(\mX_0, \Delta_{\tc, 0})$ is also a K-semistable log Fano pair. Then $\Fut(\mX, \Delta_\tc)=0$.
\end{corollary}

\begin{proof}
We only assume $(X, \Delta)$ is a log Fano $\bR$-pair since the $\bQ$-coefficients case is well known. By Proposition \ref{prop: fut=beta}, we see $\Fut(\mX, \Delta_\tc)\geq 0$.
If $\Fut(\mX, \Delta_\tc)>0$, one could derive a contradiction as in the $\bQ$-coefficients case by \cite{PT06}. We sketch it for the readers' convenience. Note that $\Fut(\mX, \Delta_\tc)$ can be reformulated as the weight of the $\bG_m$-action on the Chow-Mumford ($\bR$-) line bundle (e.g. \cite[Section 8.1]{LZ24}), where the $\bG_m$-action comes from a one parameter subgroup $\rho: \bG_m\to \Aut(\mX_0, \Delta_{\tc,0})$. Since $\rho^{-1}$ induces a product test configuration of $(\mX_0, \Delta_{\tc,0})$ with the generalized Futaki invariant being $-\Fut(\mX, \Delta_\tc)<0$, we get a contradiction to $(\mX_0, \Delta_{\tc,0})$ being K-semistable.
\end{proof}

The following result could be viewed as a converse statement (in some sense) to Proposition \ref{prop: QR}.

\begin{proposition}\label{prop: strong QR}
Let $(X, \Delta:=\sum_{j=1}^k r_j D_j)$ be a log Fano pair, where $D_j$'s are effective Weil divisors and $(r_1,...,r_k)$ is irrational. Then there exist a rational polytope $P$ containing $(r_1,...,r_k)$ as an interior point and two positive real numbers $0<\delta_0, \epsilon_0<1$ such that we have the following conclusions:
\begin{enumerate}
\item $(X, \sum_{j=1}^k x_jD_j)$ is log Fano with $ \delta(X, \sum_{j=1}^k x_jD_j)\geq \delta_0$ for any $(x_1,...,x_k)\in P$, and $\vol(-K_X-\sum_{j=1}^k x_jD_j)\geq \epsilon_0$ for any $(x_1,...,x_k)\in P$; moreover, if $\delta(X, \Delta)<1$, we may choose smaller $P$ such that $ \delta(X, \sum_{j=1}^k x_jD_j)< 1$ for any $(x_1,...,x_k)\in P$;
\item for any $(x_1,...,x_k)\in P$ and any special test configuration $(\mX, \sum_{j=1}^kx_j\mD_j)\to \bA^1$ of $(X, \sum_{j=1}^kx_jD_j)$ with $\delta(\mX_0, \sum_{j=1}^kx_j\mD_{j,0})\geq \delta_0$,  we have $(\mX, \sum_{j=1}^kr_j\mD_j)\to \bA^1$ is also a special test configuration of $(X, \Delta)$.
\end{enumerate}
\end{proposition}

\begin{proof}
The existence of the rational polytope $P$ follows from Proposition \ref{prop: QR}, and the existence of the two positive real numbers in statement (1) follows from the continuity of the $\delta$-invariant (e.g. Proposition \ref{prop: continuity}) and the continuity of the volume.

Suppose the rank of $(r_1,...,r_k)$ over $\bQ$ is $l\in \bZ_{+}$, we may assume $\dim P=l$.

We define a set of couples, denoted by $\mS$. We say a couple $(Y, \sum_{j=1}^kB_j)$ is contained in $\mS$ if and only if the following condition is satisfied: there exists some vector $(x_1,...,x_k)\in P$ such that $(Y, \sum_{j=1}^kx_j B_j)$ is the central fiber of a special test configuration of $(X, \sum_{j=1}^kx_jD_j)$ with $\delta(Y, \sum_{j=1}^kx_j B_j)\geq \delta_0$.
By \cite[Propositions 7.1 and 7.3]{LZ24}, $\mS$ is log bounded and there exists a finite rational polytope chamber decomposition of $P$, denoted by $P=\cup_s P_s$, such that for any index $s$ and any face $F$ of $P_s$, we have the following conclusions: 
\begin{enumerate}
\item for any $(Y, \sum_{j=1}^kB_j)\in \mS$, if $(Y, \sum_{j=1}^kc_jB_j)$ is log Fano for some $(c_1,...,c_k)\in F^\circ$, then $(Y, \sum_{j=1}^kx_jB_j)$ is log Fano for any $(x_1,...,x_k)\in F^\circ$;
\item if $(\mX, \sum_{j=1}^k c_j\mD_j)\to \bA^1$ is a special test configuration of $(X, \sum_{j=1}^kc_jD_j)$ for some $(c_1,...,c_k)\in F^\circ$ with $\delta(\mX_0, \sum_{j=1}^k c_j\mD_{j,0})\geq \delta_0$, then $(\mX, \sum_{j=1}^k x_j\mD_j)\to \bA^1$ is a special test configuration for any $(x_1,...,x_k)\in F^\circ$.
\end{enumerate}
Since $\rank_\bQ (r_1,...,r_k)=\dim P=l$, 
there must be some chamber $P_s$ (in the chamber decomposition of $P$) such that $P_s$ contains $(r_1,...,r_k)$ as an interior point.  Let $Q\subset P_s^\circ$ be a rational polytope which contains $(r_1,...,r_k)$ as an interior point, by the second condition listed above and replacing $P$ with $Q$, we clearly obtain the second statement in Proposition \ref{prop: strong QR}.
\end{proof}

\subsection{Complements}

Given a log pair $(X, \Delta)$ with klt singularities. We say an effective $\bR$-divisor $\Gamma$ is an \textit{$\bR$-complement }of $(X, \Delta)$ if $K_X+\Delta+\Gamma\sim_\bR 0$ and $(X, \Delta+\Gamma)$ admits log canonical singularities; if $\Delta$ is a $\bQ$-divisor and $\Gamma$ is also a $\bQ$-divisor, we say $\Gamma$ is a \textit{$\bQ$-complement}. Suppose $\Gamma$ is a $\bQ$-complement satisfying $N(K_X+\Delta+\Gamma)\sim 0$ for some positive integer $N$, we say $\Gamma$ is an $N$-complement of $(X, \Delta)$.
For an $\bR$-complement $\Gamma$, we define the \textit{dual complex} of $(X, \Delta+\Gamma)$ as follows:
$$\mathcal{DMR}(X, \Delta+\Gamma):=\{v\in \Val_X^\circ\ |\ \textit{$A_{X, \Delta+\Gamma}(v)=0$ {\rm and} $A_{X, \Delta}(v)=1$}\}. $$
Note that the space of all lc places of $(X, \Delta+\Gamma)$ is a cone over
$\mathcal{DMR}(X, \Delta+\Gamma)$, denoted by $\LCP(X, \Delta+\Gamma)$. 

We also define \textit{the complement of a Fano type variety}. Given a Fano type variety $X$, we say an effective $\bQ$-divisor $\Gamma$ is a \textit{$\bQ$-complement (resp. $N$-complement) of $X$} if $K_X+\Gamma\sim_\bQ 0$ (resp. $N(K_X+\Gamma)\sim 0$) and $(X, \Gamma)$ is log canonical.

 The following concept on special $\bQ$-complements is invented by \cite{LXZ22}, which will play an important role in Section \ref{sec: optimal deg}.

\begin{definition}[\cite{LXZ22}]\label{def: special comp}
Let $(X,\Delta)$ be a log Fano $\bQ$-pair. Let $f: (Y, E) \to (X,\Delta)$ be a log smooth model. An effective $\bQ$-divisor $\Gamma$ is called a \emph{special $\bQ$-complement} of $(X,\Delta)$ with respect to $(Y,E)$ if $(X, \Delta+\Gamma)$ is log canonical, $K_X + \Delta+\Gamma \sim_{\bQ} 0$, and $f_*^{-1}\Gamma\geq G$ for some ample effective $\bQ$-divisor $G$ on $Y$ not containing any strata of $(Y,E)$. 
\end{definition}

\subsection{Boundedness and log Boundedness}

\begin{definition}
Let $\mP$ be a set of projective varieties of dimension $d$. We say $\mP$ is bounded if there exists a projective morphism $ \mY\to T$ between schemes of finite type such that for any $X\in \mP$, there exists a closed point $t\in T$ such that $X\cong \mY_t$. Let $\mP'$ be a set of projective couples of dimension $d$. We say 
$\mP'$ is log bounded if there exist a projective morphism $\mY\to T$ between schemes of finite type and a reduced Weil divisor $\mD$ on $\mY$ such that for any $(X, \Delta)\in \mP$, there exists a closed point $t\in T$ such that $(X, \Delta)\cong (\mY_t, \mD_t)$. 
\end{definition}

We recall some results on boundedness of log Fano pairs. 

\begin{theorem}{\rm (\cite{Birkar19, Birkar21})}\label{thm: BAB}
Fix a positive integer $d$ and positive real number $\epsilon$. The following set of algebraic varieties is bounded:
$$\{X\ | \ \text{$(X, \Delta)$ is an $\epsilon$-lc log Fano pair of dimension $d$ for some $\bR$-divisor $\Delta$ on $X$}\} .$$
\end{theorem}

Based on the above boundedness, Jiang (\cite{Jiang20}) proves the boundedness of K-semistable Fano varieties with fixed dimension and volume, and this is an important step to construct the K-moduli space (e.g. \cite{Xu21}).

\begin{theorem}{\rm \cite{Jiang20}}\label{thm: bddkss}
Fix a positive integer $d$ and two positive real numbers $v$ and $a$. The following set of algebraic varieties is bounded:
$$\{X\ |\ \text{$X$ is a  Fano variety of dimension $d$ with $\vol(-K_X)\geq v$ and $\delta(X)\geq a$}\}. $$
\end{theorem}

Via the concept of normalized volume,  \cite[Theorem 6.13]{LLX19} establishes the relationship between normalized volumes and singularities, and this provides another proof of boundedness of K-semistable Fano varieties. We will recall this method below. We say $(X, \Delta; x)$ is a klt germ if $(X, \Delta)$ is klt near $x\in X$.
\begin{definition}
Let $(X, \Delta; x)$ be a klt germ of dimension $d$, the normalized volume of $(X, \Delta;x)$ is defined as follows: 
$$\widehat{\vol}(X, \Delta; x):=\inf_{v\in \Val_{X, x}^*} A_{X, \Delta}(v)^d\cdot\vol(v). $$
Here $\Val_{X, x}^*$ is the valuation space consisting of all non-trivial valuations over $X$ centered at $x$,  $A_{X, \Delta}(v)$ is the log discrepancy of $v$ with respect to $(X, \Delta)$ (e.g. \cite{JM12}), and $\vol(v)$ is defined as follows:
$$\vol(v):=\lim_{m\to \infty} \frac{\length (\mO_{X,x}/\ka_m(v))}{m^d/d!}, $$
where $\ka_m(v):=\{f\in \mO_{X, x}\ |\ \text{$v(f)\geq m$}\}$.
\end{definition}

By \cite{Blum18}, the infimum in the definition of the normalized volume can be replaced by minimum, and by \cite{Xu20, XZ21}, the minimizer is a unique quasi-monomial valuation up to scaling. 

\begin{proposition}{\rm (\cite{LLX19})}\label{prop: BL}
Fix a positive integer $d$ and a positive real number $\eta$. Then there exists a positive real number $\epsilon$ depending only on $d$ and $\eta$ such that  any klt germ $(X, \Delta; x)$ of dimension $d$ satisfying $\widehat{\vol}(X, \Delta;x)\geq \eta$ is $\epsilon$-lc in a neighborhood of $x$. 
\end{proposition}

The following result is a well-known corollary of the above proposition.

\begin{corollary}\label{cor: bdd}
Fix a positive integer $d$ and a positive real number $\eta$. Then the following set of algebraic varieties is bounded:
$$\{X \ |\  \text{$(X, \Delta)$ is a log Fano pair of dimension $d$ with $\widehat{\vol}(X, \Delta;x)\geq \eta$ for any $x\in X$}\}. $$
In particular, fixing another two positive real numbers $v$ and $a$, the following set of algebraic varieties is bounded:
$$\{X\ | \ \text{$(X, \Delta)$ is a log Fano pair of dimension $d$ with $\vol(-K_X-\Delta)\geq v$ and $\delta(X, \Delta)\geq a$}\}. $$
\end{corollary}

\begin{proof}
For the first statement, by Proposition \ref{prop: BL}, there exists a positive real number $\epsilon>0$ depending only on $d$ and $\eta$ such that $(X, \Delta)$ admits $\epsilon$-lc singularities. This implies the desired boundedness by Theorem \ref{thm: BAB}.
We turn to the second statement. Let $(X, \Delta)$ be a log Fano pair of dimension $d$ with $\vol(-K_X-\Delta)\geq v$ and $\delta(X, \Delta)\geq a$. By \cite[Proposition 3.3]{LZ24} and Proposition \ref{prop: continuity}, we can perturb the coefficients of $\Delta$ to get a $\bQ$-divisor $\Delta'$ with the same support as $\Delta$ such that $(X, \Delta')$ is a log Fano $\bQ$-pair, $\vol(-K_X-\Delta') \geq \frac{v}{2}$ and $\delta(X,\Delta') \geq \frac{a}{2}$. By \cite[Theorem D]{BJ20}, for any closed point $x\in X$ we have
$$\widehat{\vol}(X, \Delta'; x)\geq \frac{d^d}{(1+d)^d}\cdot \vol(-K_X-\Delta')\cdot \delta(X, \Delta')^d\geq  \frac{d^d\cdot a^d\cdot v}{(1+d)^d\cdot 2^{d+1}}.$$
Thus the boundedness follows from the first statement.
\end{proof}

\section{Special test configuration theory}\label{sec: LX}

Given a log Fano $\bQ$-pair $(X, \Delta)$, the special test configuration theory is well established in \cite{LX14}, i.e. to test K-semistability of $(X, \Delta)$, it suffices to confirm the non-negativity of the generalized Futaki invariants of all special test configurations of $(X, \Delta)$. In this section, we aim to prove the corresponding result for log Fano $\bR$-pairs. However, we do not plan to rebuild the whole program in \cite{LX14} since we take Fujita-Li's criterion as the definition for K-semistability (e.g. Definition \ref{def: Fujita-Li}), rather than by test configurations. 

\begin{definition}\label{def: optimal deg}
Given a log Fano $\bQ$-pair $(X, \Delta)$ with $\delta(X, \Delta)\leq 1$. Let $E$ be a prime divisor over $X$ computing $\delta(X, \Delta)$ (such $E$ always exists due to \cite{LXZ22}), then $E$ induces a special test configuration $(\mX, \Delta_\tc)\to \bA^1$ with $\delta(\mX_0, \Delta_{\tc, 0})=\delta(X, \Delta)$ (e.g. \cite{BLZ22}). We call $(\mX, \Delta_\tc)\to \bA^1$ or $(\mX_0, \Delta_{\tc, 0})$ an \textit{optimal destabilization} of $(X, \Delta)$.
\end{definition}

\begin{theorem}\label{thm: LX}
Let $(X, \Delta:=\sum_{j=1}^k r_j D_j)$ be a log Fano pair of dimension $d$, where $D_j$'s are effective Weil divisors and $(r_1,...,r_k)$ is irrational. Then $(X, \Delta)$ is K-semistable if and only if $\Fut(\mX, \Delta_{\tc})\geq 0$ for any special test configuration $(\mX, \Delta_{\tc})\to \bA^1$ of $(X, \Delta)$.
\end{theorem}

\begin{proof}
Suppose $(X, \Delta)$ is K-semistable. By Proposition \ref{prop: fut=beta}, we see that $\Fut(\mX, \Delta_{\tc})\geq 0$ for any special test configuration $(\mX, \Delta_{\tc})\to \bA^1$ of $(X, \Delta)$.

Conversely, suppose $\Fut(\mX, \Delta_{\tc})\geq 0$ for any special test configuration $(\mX, \Delta_{\tc})\to \bA^1$ of $(X, \Delta)$. We aim to show that $(X, \Delta)$ is K-semistable, i.e. $\delta(X, \Delta)\geq 1$.
Assume to the contrary that $\delta(X, \Delta)<1$. We choose a rational polytope $P$ and two positive real numbers $\delta_0, \epsilon_0$ as in Proposition \ref{prop: strong QR}. Let $\{(r_{i1},...,r_{ik})\}_{i=1}^\infty\subset P$ be a sequence of rational vectors tending to $(r_1,...,r_k)$. 
Denote by $\Delta_i:=\sum_{j=1}^k r_{ij}D_j$. By the continuity of the $\delta$-invariant (Proposition \ref{prop: continuity}), we have
$$\delta(X, \Delta)=\lim_{i\to\infty} \delta(X, \Delta_i)<1.$$ 
We may take $\delta_0 =\frac{\delta(X,\Delta)}{2}$.
Thus there exists $i_0\in \bN$ such that $\delta_0<\delta(X,\Delta_i)<1$ for every $i \geq i_0$. 
For each index $i \geq i_0$, let $E_i$ be a prime divisor over $X$ which induces an optimal destabilization of $(X, \Delta_i)$, denoted by $(\mX^{(i)}, \sum_{j=1}^k r_{ij}\mD_j^{(i)})\to \bA^1$. Denote by $\Phi:=\{E_i\}_{i\geq i_0}$. It is clear to see 
$$\delta(X, \Delta_i)=\frac{A_{X, \Delta_i}(E_i)}{S_{X, \Delta_i}(E_i)}=\inf_{E\in \Phi} \frac{A_{X, \Delta_i}(E)}{S_{X, \Delta_i}(E)}.$$
Thus we have the following $(\star)$:
$$\delta(X, \Delta)=\lim_{i\to\infty} \delta(X, \Delta_i)=\lim_{i\to\infty} \inf_{E\in \Phi}\frac{A_{X, \Delta_i}(E)}{S_{X, \Delta_i}(E)}. $$
We claim that it suffices to show the following equality $(\star\star)$:
$$\lim_{i\to\infty} \inf_{E\in \Phi}\frac{A_{X, \Delta_i}(E)}{S_{X, \Delta_i}(E)}= \inf_{E\in \Phi} \frac{A_{X, \Delta}(E)}{S_{X, \Delta}(E)}.$$
Granting $(\star\star)$, we see the following formula by combining $(\star)$ with $(\star\star)$:
$$\delta(X, \Delta)=\inf_{E\in \Phi} \frac{A_{X, \Delta}(E)}{S_{X, \Delta}(E)}.$$
Since $\delta(\mX_0^{(i)}, \sum_{j=1}^k r_{ij}\mD_{j,0}^{(i)}) = \delta(X, \Delta_i) >\delta_0$ for every $i\geq i_0$, by Proposition \ref{prop: strong QR} we have that  $(\mX^{(i)}, \sum_{j=1}^k r_{j}\mD_j^{(i)})\to \bA^1$ is a special test configuration of $(X, \Delta)$, i.e. every  $E\in \Phi$ is also a special divisor for $(X, \Delta)$. Applying Proposition \ref{prop: fut=beta}, we see $A_{X, \Delta}(E)-S_{X, \Delta}(E)\geq 0$ for any $E\in \Phi$. Thus 
$\delta(X, \Delta)\geq 1$, leading to a contradiction.

Finally, the equality $(\star\star)$ is implied by the same argument in Proposition \ref{prop: continuity}.
The proof is complete.
\end{proof}

The existence of a prime divisor computing the $\delta$-invariant of a log Fano $\bQ$-pair (as in Definition \ref{def: optimal deg}) 
plays an essential role in the proof of Theorem \ref{thm: LX}. In the next section, we will show the existence of a prime divisor computing the $\delta$-invariant of a log Fano $\bR$-pair. This will lead to another proof of Theorem \ref{thm: LX} without using the approximation by $\bQ$-coefficients (see Section \ref{subsec: fg}).
The following result corresponds to the last step in \cite{LX14} of decreasing the generalized Futaki invariants, which is already enough for  later applications (in Section \ref{sec: Theta} and Section \ref{sec: S}).

\begin{proposition}\label{prop: weak LX}
Let $(X, \Delta)$ be a log Fano $\bR$-pair of dimension $d$. 
Suppose $(X, \Delta)$ is K-semistable and $(\mX, \Delta_\tc)\to \bA^1$ is a weakly special test configuration of $(X, \Delta)$ with $\Fut(\mX, \Delta_\tc)=0$. Then the test configuration is special.
\end{proposition}

\begin{proof}
By the same proof of \cite[Proposition 2.49]{Xu24}, up to a finite base change and a standard tie-breaking argument, we could cook up a special test configuration $(\mX^s, \Delta_\tc^s)\to \bA^1$ of $(X, \Delta)$ such that $A_{\mX, \Delta_\tc}(\mX_0^s)=1$. By the same proof in step 3 of \cite[Proof of Theorem 2.51]{Xu24}, we see $\Fut(\mX^s, \Delta_\tc^s)\leq \Fut(\mX, \Delta_\tc)$ up to a positive multiple, and the equality holds if and only if $(\mX, \Delta_\tc)$ is itself a special test configuration.
The proof is complete.
\end{proof}

\begin{remark}\label{rem: weak LX}
Suppose $(X, \Delta)$ is a K-semistable log Fano $\bR$-pair. From the proof of Proposition \ref{prop: weak LX}, we see that $\Fut(\mX, \Delta_\tc)\geq 0$ for any weakly special test configuration $(\mX, \Delta_\tc)$ of $(X, \Delta)$. However, here we do not attempt to prove that $\Fut(\mX, \Delta_\tc)\geq 0$ for any test configuration $(\mX, \Delta_\tc)$ of $(X, \Delta)$ under the condition that $(X, \Delta)$ is K-semistable (though this is expected to be true).
The reason is we obtain Theorem \ref{thm: LX} by a way avoiding developing the whole program in \cite{LX14}. We leave this to the readers who are interested.

\end{remark}

\section{Optimal destabilization}\label{sec: optimal deg}

In this section, we will establish the optimal destabilization theory for a log Fano $\bR$-pair. More concretely, for a given log Fano $\bR$-pair $(X, \Delta)$ with $\delta(X, \Delta)\leq 1$, we will show there exists a divisorial valuation computing $\delta(X, \Delta)$, and any such divisorial valuation is special and induces a special test configuration keeping the $\delta$-invariant on the central fiber.

\subsection{Finite generation}\label{subsec: fg}

In this subsection, we show the existence of a prime divisor computing the $\delta$-invariant, and show that the divisor is special. We will take a similar strategy as in \cite{LXZ22} and a priori prove the finite generation property. Before this, we first make some preparations. All the results in this subsection are well-known for log Fano $\bQ$-pairs by \cite{LXZ22}, so we only consider log Fano $\bR$-pairs in this subsection.

\begin{proposition}\label{prop: qm}
Let $(X, \Delta)$ be a log Fano $\bR$-pair of dimension $d$ with $\delta(X, \Delta)< \frac{d+1}{d}$. Then there exists a positive integer $N$ depending only on $d$ such that 
$$\delta(X, \Delta)=\inf_{E}\frac{A_{X, \Delta}(E)}{S_{X, \Delta}(E)}, $$
where $E$ runs over all lc places of $N$-complements of $X$; moreover, there exists a quasi-monomial valuation $v$ computing $\delta(X, \Delta)$.
\end{proposition}

\begin{proof}
The strategy is similar to the proof of \cite[Theorem 4.6]{BLX22}. For each large positive integer, there exists a prime divisor $E_m$ and an $m$-basis type divisor $B_m$ such that $S_m(E_m)=\ord_{E_m}(B_m)$ and $E_m$ computes $\delta_m(X, \Delta)$, i.e. 
$$\delta_m:=\delta_m(X, \Delta)=\lct(X, \Delta; B_m)=\frac{A_{X, \Delta}(E_m)}{S_m(E_m)}.$$
Taking a sufficiently general $D\in |-K_X-\Delta|_\bR$,  we may choose $B_m$ satisfying $\ord_{D}(B_m)=S_m(D)$ (e.g. Proposition \ref{prop: compatible}). Denote by $B_m':=B_m-S_m(D)D$. Since $\lim_{m\to \infty}S_m(D)=\frac{1}{d+1}$, we see that $E_m$ is an lc place of $(X, \Delta+\delta_m B_m')$ and 
$$\delta_m B_m'\sim_\bR \delta_m (1-S_m(D))(-K_X-\Delta),$$
where $\lim_{m\to \infty} \delta_m(1-S_m(D)) = \delta(X,\Delta) (1-\frac{1}{d+1})<1$. Thus for $m\gg 1$ there exists an $\bR$-complement $\Gamma_m$ of $(X, \Delta)$ such that $E_m$ is an lc place of $(X, \Delta+\Gamma_m)$. By the same argument  in the proof of \cite[Proposition 5.1]{LZ24}, $E_m$ is an lc place of some $N$-complement of the Fano type variety $X$, where $N$ depends only on the dimension $d$. In summary, we have $\delta(X, \Delta)=\inf_E \frac{A_{X, \Delta}(E)}{S_{X, \Delta}(E)}$, where $E$ runs over all lc places of $N$-complements of $X$.

The rest follows the proof of \cite[Theorem 4.6]{BLX22} by replacing \cite[Proposition 4.2]{BLX22} with \cite[Proposition 5.2]{LZ24}.
\end{proof}

\begin{proposition}\label{prop: special R}
Let $(X, \Delta)$ be a log Fano $\bR$-pair of dimension $d$ with $\delta(X,\Delta)<\frac{d+1}{d}$. Let $v$ be a quasi-monomial valuation computing $\delta:=\delta(X, \Delta)$ and let $\alpha \in (0, \min \{\frac{\delta}{d+1}, 1-\frac{d\delta}{d+1}\})$. 
Then for any $B\in |-K_X-\Delta|_{\bR}$, there exists an $\bR$-complement $\Gamma$ of $(X, \Delta+\alpha B)$ such that $v$ is an lc place of $(X, \Delta+\alpha B+\Gamma)$. 
\end{proposition}

\begin{proof}
We follow the same proof as \cite[Lemma 3.1]{LXZ22}.
After a rescaling, we may assume $A_{X,\Delta}(v) = 1$. We pick a log smooth model $\pi: (Y, E)\to (X,\Delta)$ such that $v\in \QM_\eta(Y,E)$, where $\eta$ is the generic point of a stratum of $E$. 

By \cite[Lemma 2.7]{LX18}, for each $\epsilon>0$, there exist divisorial valuations $v_i\in \QM_{\eta}(Y,E)$ and $q_1,\cdots, q_r\in \bZ_{>0}$ such that $v$ is in the convex cone generated by $v_i$, each $v_i = q_i^{-1} \ord_{F_i}$ for some divisor $F_i$ over $X$, and $|v_i - v|<q_i^{-1} \epsilon$ for every $i$. Let $\ka_{\bullet}=\ka_{\bullet}(v)$ be the ideal sequence of $v$. By the same argument as in \cite[Proof of Lemma 3.1]{LXZ22}, there exists $C>0$ and $\epsilon_1>0$ such that whenever $\epsilon \in (0, \epsilon_1)$ there exists $\epsilon_0\in (0,1)$ depending on $\epsilon$ such that
\[
A_{X,\Delta+\ka_{\bullet}^{1-\epsilon_0}}(F_i) < 2C\epsilon.
\]
Let $B'\in |-K_X-\Delta|_{\bR}$ be general so that $\Supp(B')$ does not contain $c_X(v)$. Let $G:=\beta B' + (1-\beta) B$ with $\beta:=\max \{0, \frac{(d+1)(\delta-1)}{\delta}\}$. Let $B_m$ be an $m$-basis type divisor of $(X, \Delta)$ compatible with both $v$ and $G$ (e.g. Proposition \ref{prop: compatible}). Thus $B_m\geq S_m(G)G$. Note that $\lim_{m\to\infty} S_m(G) = \frac{1}{d+1}$ (e.g. Section \ref{subsec: S-inv}).
Denote by $B_m':=B_m - S_m(G)\beta B'$.  We may choose a sequence of real numbers $\delta_m\in (0, \delta_m(X,\Delta))$ such that $\lim_{m\to\infty} \delta_m = \delta$ and $\delta_m(1-\beta S_m(G))<1$ for every $m\gg 1$. Thus we have 
$$B_m'\sim_{\bR} (1-\beta S_m(G)) (-K_X-\Delta)\quad \text{and} \quad\delta_m B_m'\geq \alpha B$$ 
for $m\gg 1$, and $(X, \Delta+ \delta_m B_m')$ is a log Fano $\bR$-pair.
Since $v$ computes $\delta(X, \Delta)$ and $B'$ does not contain $c_X(v)$, we see the following for $m\gg 1$:
$$\delta_mv(B_m')=\delta_mv(B_m)\geq (1-\epsilon_0)\delta(X, \Delta)S_{X, \Delta}(v)=(1-\epsilon_0)A_{X, \Delta}(v)=1-\epsilon_0. $$
Thus for any $\epsilon\in (0, \min\{\frac{1}{2C},\epsilon_1\})$, we have
\[
A_{X, \Delta+ \delta_m B_m'}(F_i) < 2C\epsilon <1.
\]
By \cite[Corollary 1.4.3]{BCHM10}, there exists a $\bQ$-factorial birational model $\tX\to X$ that precisely extracts all the $F_i$'s, and we have $(\tX, \tilde{\Delta} + \delta_m \tB_m' + \sum_{i=1}^{r} (1-a_i)F_i)$ is a weak log Fano $\bR$-pair as the crepant pull back of $(X, \Delta+\delta_m B_m')$, where $a_i = A_{X, \Delta+ \delta_m B_m'}(F_i) < 2C\epsilon$. Thus we know that $(\tX, \tilde{\Delta} + \delta_m \tB_m' + \sum_{i=1}^{r} (1-a_i)F_i)$ has an $\bR$-complement, hence so does $(\tX, \tilde{\Delta} + \alpha \tB +\sum_{i=1}^{r} (1-a_i)F_i)$. By the $\bR$-coefficients generalization of \cite[Lemma 3.2]{LXZ22} (note that it can be generalized to $\bR$-coefficients as both the MMP for Fano type varieties and the ACC for LCT work for $\bR$-pairs according to \cite{BCHM10, HMX14}), there exists $\epsilon_2\in (0,1)$ depending only on $d,~\alpha$ and the coefficients of $\Delta$ and $B$, such that for any $\epsilon \in (0, \frac{\epsilon_2}{2C})$, the pair $(\tX, \tilde{\Delta} + \alpha \tB + \sum_{i=1}^r F_i)$ admits an $\bR$-complement as each $a_i <\epsilon_2$. By the push-forward to $X$, we have an $\bR$-complement $\Gamma$ of $(X, \Delta+\alpha B)$ such that  each $F_i$ is an lc place of $(X, \Delta+\alpha B+\Gamma)$. Since the log discrepancy function $A_{X, \Delta+\alpha B+\Gamma}(\cdot)$ is convex on $\QM_{\eta}(Y,E)$, we know that $v$ must be an lc place of $(X, \Delta+\alpha B+\Gamma)$ as it is a convex combination of $F_i$'s. 
\end{proof}

\begin{proposition}\label{prop: special Q}
Let $(X, \Delta)$ be a log Fano $\bR$-pair of dimension $d$ with $\delta(X, \Delta)< \frac{d+1}{d}$. Let $v$ be a quasi-monomial valuation computing $\delta(X, \Delta)$.
Then there exist $\epsilon_0>0$, a $\bQ$-divisor $\Delta'\geq (1+\epsilon_0)\Delta$ with $(X,\Delta')$ a log Fano $\bQ$-pair, a log smooth model $(Y, E) \to (X, \Delta')$ and a special $\bQ$-complement $\Gamma''$ of $(X,\Delta')$ with respect to $(Y,E)$ such that $v\in \QM(Y,E) \cap \LCP(X, \Delta'+\Gamma'')$. 
\end{proposition}

\begin{proof}
By the proof of Proposition \ref{prop: special R}, there exists a log smooth model $\pi: (Y, E)\to (X, \Delta)$ such that $v\in \QM_\eta(Y, E)$, where $\eta$ is the generic point of a stratum of $E$; moreover, there exists a real number $\alpha>0$ such that for any $B\in |-K_X-\Delta|_\bR$, $(X, \Delta+\alpha B)$ admits an $\bR$-complement $\Gamma$  such that $v$ is contained in a simplicial cone of $\QM_\eta(Y, E)\cap \LCP(X, \Delta+\alpha B+\Gamma)$ generated by prime divisors $F_1,...,F_r$ over $X$. Let $L$ be a very ample line bundle on $Y$ and $G\in |L|$ be a general element whose support does not contain any strata of $E$. Since $-K_X-\Delta$ is ample, we may take $B\in |-K_X-\Delta|_\bR$ to satisfy $B>b(\pi_*G+\Delta)$ for some $b>0$. One can easily choose a $\bQ$-divisor $\Delta'$ satisfying $\Delta+\frac{\alpha}{4}B\leq \Delta'\leq\Delta+\frac{\alpha}{2} B$ and $(X, \Delta')$ is a log Fano $\bQ$-pair. Moreover, we have
\begin{align*}
\Delta+\alpha B+\Gamma = \Delta'+(\Delta+\frac{\alpha}{2}B-\Delta')+\frac{\alpha}{2} B +\Gamma.
\end{align*}
Let $\epsilon_0:=\frac{\alpha b}{4}$ and denote $\Gamma':=(\Delta+\frac{\alpha}{2}B-\Delta')+\frac{\alpha}{2} B +\Gamma$. Then we see
$$\Delta+\alpha B+\Gamma=\Delta'+\Gamma',\ \ \ \Delta'\geq (1+\epsilon_0)\Delta, \quad \text{and}\quad \Gamma'\geq \frac{\alpha b}{2}\pi_*G.$$ 
By \cite[Corollary 1.4.3]{BCHM10}, there exists a $\bQ$-factorial birational model $f: \tilde{X}\to (X, \Delta'+\Gamma')$ which precisely extracts all $F_i$'s, and we write
$$K_{\tX}+\widetilde{\Delta'}+\widetilde{\Gamma'}+\sum_{i=1}^r F_j=f^*(K_X+\Delta'+\Gamma'). $$
Here we write $\widetilde{\Delta'},\widetilde{\Gamma'}, \widetilde{\pi_*G}$ for the birational transform of $\Delta', \Gamma', \pi_*G$ on $\tX$. Note that $\tX$ is Fano type. Let $0<b_0<\frac{\alpha b}{2}$ be a rational number, then we run MMP on 
$$-K_{\tX}-\widetilde{\Delta'}-b_0\widetilde{\pi_*G}-\sum_{i=1}^r F_i\sim_\bR \widetilde{\Gamma'}-b_0\widetilde{\pi_*G} \geq 0$$
which terminates with a good minimal model $\phi: \tX\dashrightarrow X'$. Note that the left hand side is a $\bQ$-divisor. By \cite[Section 6.1]{Birkar19}, $(\tX, \widetilde{\Delta'}+b_0\widetilde{\pi_*G}+\sum_{i=1}^r F_i)$ admits a $\bQ$-complement. By the push-forward to $X$, one clearly gets a desired special $\bQ$-complement $\Gamma''$ for $(X, \Delta')$.
\end{proof}

We are ready to prove the finite generation property. We first fix some notation. Let $X$ be a projective normal variety and $D_j$'s for $1\leq j\leq k$ are effective Weil divisors on $X$, we denote $\Delta(\fa):=\sum_{j=1}^ka_j D_j$ for $\fa=(a_1,...,a_k)\in \bR^k_{\geq 0}$. Suppose $(X, \Delta(\fa'))$ is a log Fano $\bQ$-pair for some $\fa':=(a_1',...,a_k')\in \bQ^k_{\geq 0}$, we write 
$$R(\fa'):=R(X, -K_X-\Delta(\fa')):=\oplus_{m\in M(\fa')} H^0(X, -m(K_X+\Delta(\fa'))),$$
where $M(\fa'):=\{m\in \bN\ |\ \textit{$-m(K_X+\Delta(\fa'))$ {\rm is Cartier}}\}$.
For a valuation $v\in \Val_X^\circ$, the graded ring associated to $v$ is denoted by 
$$\gr_v R(\fa'):=\oplus_{m\in M(\fa')}\oplus_{\lambda\in \bR_{\geq 0}} \ka_{m,\lambda}(v)/\ka_{m, >\lambda}(v),$$
where $\ka_{m,\lambda}(v):=\{s\in R(\fa')_m\ |\ \textit{$v(s)\geq \lambda$}\}$. Let $Y\to X$ be a proper normal birational model, we denote $c_Y(v)$ to be the center of $v$ on $Y$.

\begin{theorem}\label{thm: fg}
Let $(X, \Delta=\Delta(\fa):=\sum_{j=1}^ka_jD_j)$ be a log Fano $\bR$-pair of dimension $d$, where $D_j$'s are effective Weil divisors and $\fa=(a_1,...,a_k)\in (0,1)^k$ is irrational. Suppose $\delta(X, \Delta)< \frac{d+1}{d}$ and let $v$ be a quasi-monomial valuation computing $\delta(X, \Delta)$. Then there exist an open neighborhood $U$ of $\fa$ contained in the log Fano polytope and a smooth cone $\Sigma$ in $\QM_{\eta}(Y,E)$ containing $v$ (as an interior point) for some log smooth model $(Y,E) \to (X,\Delta)$ with $c_Y(v)=\eta$, such that for every $\fa'\in U\cap \bQ^k$ and every $v_w\in \Int(\Sigma)$ with weight $w\in \bR_{>0}^k$, the graded ring $\gr_{v_w} R(\fa')$ is finitely generated and independent of the choice of $w$.  
\end{theorem}

\begin{proof}
By Proposition \ref{prop: special Q}, there exist $\epsilon_0>0$, a $\bQ$-divisor $\Delta'\geq (1+\epsilon_0)\Delta$, a log smooth model $(Y,E)\to (X,\Delta')$, and a special $\bQ$-complement $\Gamma$ of $(X,\Delta')$ with respect to $(Y,E)$ such that $v\in \QM_{\eta}(Y,E)\cap \mathrm{LCP}(X,\Delta'+\Gamma)$, where $\eta := c_Y(v)$. Since $\Delta'\geq (1+\epsilon_0)\Delta$, there exists a neighborhood $U$ of $\fa$ in the log Fano polytope such that $\Delta'\geq \Delta(\fa')$ for any $\fa'\in U$ (e.g. Proposition \ref{prop: QR}). 

Let $F_1,\cdots, F_r\in \QM_{\eta}(Y,E)\cap \mathrm{LCP}(X,\Delta'+\Gamma)$ be prime divisors over $X$ so that they generate a smooth cone $\Sigma$ in $\QM_{\eta}(Y,E)$ containing $v$. Then we know that $\Gamma+(\Delta'-\Delta(\fa'))$ is a special $\bQ$-complement of $(X,\Delta(\fa'))$ for any $\fa'\in U\cap \bQ^k$. 
By \cite[Theorem 5.26]{Xu24}, for every $\fa'\in U\cap \bQ^k$ there exists a dlt Fano type model $(Y',E')\to (X, \Delta(\fa'))$ such that $E':=\sum_{i=1}^rF_i'$ is the sum of the birational transform of $F_1,\cdots, F_r$, and the toroidal structure of $(Y',E')$ at $\cap_{i=1}^r F_i'$ coincides with the one from $(Y,E)$. Note that $(Y',E')$ may depend on the choice of $\fa'$. Thus by \cite[Theorem 5.20]{Xu24} we know that $\gr_v R(\fa')$ is finitely generated, and $\gr_v R(\fa')\cong \gr_{v_w} R(\fa')$ for any $w\in \bR_{>0}^r$. 
The proof is finished. 
\end{proof}

As a corollary, we obtain the following existence of a special divisorial minimizer for the $\delta$-invariant.

\begin{theorem}\label{thm: special minimizer}
Let $(X, \Delta)$ be a log Fano $\bR$-pair of dimension $d$ with $\delta(X, \Delta)< \frac{d+1}{d}$. Then there exists a prime divisor over $X$ computing $\delta(X, \Delta)$. Moreover, any such prime divisor is special.
\end{theorem}

\begin{proof}
We write $\Delta=\Delta(\fa)=\sum_{j=1}^ka_jD_j$ as in Theorem \ref{thm: fg}.
Let $v_0$ be a quasi-monomial valuation computing $\delta(X,\Delta)$ whose existence follows from Proposition \ref{prop: qm}. We follow the notation in Theorem \ref{thm: fg}. By \cite[Theorem 4.5]{LXZ22},  in the small cone $\Sigma$ (containing $v_0$ as an interior point) of $\QM_{\eta}(Y,E)$, the function
$S_{X, \Delta(\fa')}(\cdot)$
is linear on $\Sigma$ for any $\fa'\in U\cap \bQ^k$. Since the $S$-invariant is continuous in the coefficients, we know that the function $S_{X,\Delta}(\cdot)$ is also linear on $\Sigma$. Since $(Y,E)\to (X,\Delta)$ is a log smooth model, we know that $A_{X,\Delta}(\cdot)$ is linear on $\Sigma$. Note that $v_0$ computing $\delta(X,\Delta)$ implies that $A_{X,\Delta}(v_0)- \delta(X,\Delta) S_{X,\Delta} (v_0) =0$. Since the function
\[
u\mapsto A_{X,\Delta}(u) - \delta(X,\Delta) S_{X,\Delta}(u)
\]
is non-negative and linear on $\Sigma$, we conclude that it must vanish on $\Sigma$. Thus any rational point in $\Sigma$ corresponds to a divisorial valuation that also computes $\delta(X,\Delta)$. Denote $\Sigma(\bQ)$ to be the set of rational points in $\Sigma$. Let $v\in \Sigma(\bQ)$ be any rational point, it suffices to show that $v$ is a special divisorial valuation for $(X, \Delta=\Delta(\fa))$. By the proof of Proposition \ref{prop: special Q}, we see that $v$ is an lc place of a special $\bQ$-complement of $(X, \Delta(\fa'))$ for any $\fa':=(a_1',...,a_k')\in U\cap \bQ^k$ with respect to $(Y, E)$. By \cite[Theorem 4.2]{LXZ22}, we know that $v$ induces a special test configuration of $(X,\Delta(\fa')=\sum_{j=1}^ka_j'D_j)$, denoted by
$$(\mX', \sum_{j=1}^ka_j'\mD_j')\to \bA^1.$$ 
Take a $\bC^*$-equivariant small $\bQ$-factorization $(\mX'', \sum_{j=1}^k\mD_j'')\to (\mX', \sum_{j=1}^k\mD_j')$  and run a $\bG_m$-equivariant MMP/$\bA^1$ for $-K_{\mX''}- \sum_{j=1}^ka_j\mD_j''$ which terminates with a canonical model $(\mX, \sum_{j=1}^ka_j\mD_j)$. 
Note that we can run this MMP as $(\mX', \sum_{j=1}^ka_j'\mD_j')$ is a special test configuration of a log Fano $\bQ$-pair $(X, \Delta(\fa'))$, hence (a small $\bQ$-factorialization of) $\mX$ is a Mori dream space over $\bA^1$. 
We claim that $(\mX, \sum_{j=1}^ka_j\mD_j)\to \bA^1$ is a family of log Fano $\bR$-pairs, which implies that it is a special test configuration of $(X, \Delta(\fa))$. 
From the MMP we know that $-K_{\mX}-\sum_{j=1}^ka_j\mD_j$ is ample. Let $\Pi$ be the polyhedron of $\bR_{\geq 0}^k$ such that $K_{\mX}+\sum_{j=1}^kx_j\mD_j$ is $\bQ$-Cartier for every $(x_1,...,x_k)\in \Pi$. Then we know that $\Pi$ is rational (e.g. \cite[Proposition 3.2]{LZ24}) and it contains $\fa=(a_1,...,a_k)$ as an interior point.  Let $V$ be an open  neighborhood of  $\fa$ in $U\cap \Pi$ such that $-K_{\mX}-\sum_{j=1}^ka_j''\mD_j$ is $\bQ$-Cartier and ample for any $\fa'':=(a_1'',...,a_k'')\in V\cap \bQ^k$. Recall that $v$ induces a special test configuration of  $(X, \Delta(\fa''))$ for any $\fa'':=(a_1'',...,a_k'')\in V\cap \bQ^k$,
which has to be $(\mX, \sum_{j=1}^ka_j''\mD_j)$ by the ampleness of  $-K_{\mX}-\sum_{j=1}^ka_j''\mD_j$. 
Thus we know that $(\mX_0, \sum_{j=1}^ka_j''\mD_{j,0})$ is a log Fano pair for any $\fa''\in V\cap \bQ^k$. 
Since $\fa$ lies in the convex hull of some finite subset of $V\cap \bQ^k$, we conclude that $(\mX_0, \sum_{j=1}^ka_j\mD_{j,0})$ is a klt log Fano $\bR$-pair. The proof is finished. 
\end{proof}

\begin{proof}[The second proof of Theorem \ref{thm: LX}]
We only show that $(X, \Delta)$ is K-semistable under the condition that 
$\Fut(\mX, \Delta_{\tc})\geq 0$ for any special test configuration $(\mX, \Delta_{\tc})\to \bA^1$ of $(X, \Delta)$. Suppose $\delta(X, \Delta)<1$,  by Theorem \ref{thm: special minimizer}, there exists a prime divisor $E$ over $X$ computing $\delta(X, \Delta)$, and $E$ induces a special test configuration of $(X, \Delta)$, denoted by  $(\mX, \Delta_\tc)\to \bA^1$. Applying Proposition \ref{prop: fut=beta}, we know $\Fut(\mX, \Delta_\tc)$ is equal to $A_{X, \Delta}(E)-S_{X, \Delta}(E)$ up to a positive multiple. This implies $\Fut(\mX, \Delta_\tc)<0$ and leads to a contradiction.
\end{proof}

\subsection{$\bT$-equivariant K-stability}

In this subsection, we aim to obtain an equivariant version of Theorem \ref{thm: special minimizer}. As an application, we will establish the equivariant analogue of Theorem \ref{thm: LX} (see Theorem \ref{thm: G-LX}).  Throughout the subsection, we fix $(X, \Delta)$ to be a log Fano pair of dimension $d$ with a torus $\bT:=\bG_m^r$-action for some $r\in \bN$. Denote $\Val_X^{\bT}$ (resp. $\DivVal_X^{\bT}$) to be the subset of $\Val_X^\circ$ consisting of $\bT$-invariant valuations (resp. $\bT$-invariant divisorial valuations). For any prime divisor $F$ over $X$ and each large positive integer $m$ we define
$$S_m^{\bT}(F):=\sup_{\text{$\bT$-invariant $B_m$}}\ord_F (B_m)\quad \text{and}\quad\delta_m^{\bT}(X, \Delta):=\inf_{\text{$\bT$-invariant $B_m$}}\lct(X, \Delta; B_m),$$
where $B_m$ runs over all $\bT$-invariant $m$-basis type divisors.

\begin{proposition}\label{prop: T-delta}
Notation as above, we have the following equalities:
$$\delta_m(X, \Delta)=\delta_m^\bT(X, \Delta)=\inf_{v\in \Val_X^{\bT}}\frac{A_{X, \Delta}(v)}{S_m^{\bT}(v)}=\inf_{E\in \DivVal_X^{\bT}}\frac{A_{X, \Delta}(E)}{S_m^{\bT}(E)}.$$
\end{proposition}

\begin{proof}
The second and third equalities follow from the same argument as \cite{BJ20, AZ22} together with the fact that the $\lct$ of a $\bT$-equivariant $\bR$-Cartier $\bR$-divisor is always computed by $\bT$-equivariant divisorial valuation, which is a consequence of the existence of $\bT$-equivariant log resolutions. Since it is clear that $\delta_m(X, \Delta) \leq \delta_m^{\bT}(X, \Delta)$, it suffices to show the reverse inequality.
Suppose that $\delta_m(X, \Delta)$ is computed by a (not necessarily $\bT$-equivariant) prime divisor $E$ over $X$. Then for any basis type divisor $B_m$ compatible with $E$ (i.e. $S_m(E)=\ord_E(B_m)$), we have 
\[
\delta_m(X, \Delta)\leq \lct(X, \Delta;B_m)\leq \frac{A_{X, \Delta}(E)}{\ord_E(B_m)} = \frac{A_{X, \Delta}(E)}{S_m(E)} = \delta_m(X, \Delta),
\]
which implies $\lct(X,\Delta;B_m) = \delta_m(X, \Delta)$.  Let $\mF$ be the filtration on $H^0(X, -m(K_X+\Delta))$ induced by $\ord_E$ (see \cite[Section 2.3 and Section 3.1]{BJ20}), and $\mF_0$ the initial degeneration of $\mF$. Then there exists a flat $\bT$-equivariant degeneration of $\mF$ to $\mF_0$. As a result, there exists a basis type divisor $B_m$ compatible with $\mF$ (i.e. compatible with $E$) that  isotrivially degenerates to a $\bT$-equivariant basis type divisor $B_{m,0}$ compatible with $\mF_0$ (see \cite[Definition 2.7 and Lemma 3.1]{AZ22}). By lower semicontinuity of $\lct$, we have 
\[
\delta_m(X, \Delta) = \lct(X,\Delta;B_m) \geq \lct(X, \Delta;B_{m,0}) \geq \delta_m^{\bT}(X, \Delta).
\]
The proof is finished.
\end{proof}

\begin{theorem}\label{thm: G-minimizer}
Let $(X, \Delta)$ be a log Fano pair which admits a $\bT$-action. Suppose $\delta(X, \Delta)< 1$. Then there exists a $\bT$-equivariant special divisorial valuation computing $\delta(X, \Delta)$.
\end{theorem}

\begin{proof}
By the proof of Proposition \ref{prop: T-delta}, for each large positive integer $m$, there exist a  $\bT$-invariant prime divisor  $E_m$ over $X$ and a $\bT$-invariant $m$-basis type divisor $B_m$ such that $E_m$ computes $\delta_m:=\delta_m(X, \Delta)=\lct(X, \Delta; B_m)<1$. 
Following the same argument as in the proof of Proposition \ref{prop: qm}, we conclude that $E_m$ is an lc place of some $N$-complement $\Gamma_m$ of $X$, where $N$ depends only on the dimension $d$. 
By the proof of \cite[Theorem 3.24]{LW24}, after possibly enlarging $N$ we may assume that $\Gamma_m$ is $\bT$-invariant. Thus the same argument as in the proof of Proposition \ref{prop: qm} implies that there exists a $\bT$-invariant quasi-monomial valuation $v$ computing $\delta(X,\Delta)$. By the proof of Proposition \ref{prop: special R} and Theorem \ref{thm: fg}, there exists a $\bT$-equivariant log smooth model $(Y,E)\to (X,\Delta)$ such that $v\in \QM_{\eta}(Y,E)$ and that $S_{X,\Delta}(\cdot)$ is linear in a simplicial cone in $\QM_{\eta}(Y,E)$ which contains $v$. Thus the proof of Theorem \ref{thm: special minimizer} implies that we can perturb $v$ to a $\bT$-invariant divisorial valuation which computes $\delta(X,\Delta)$ as well. The proof is finished.
\end{proof}

\begin{theorem}\label{thm: G-LX}
Let $(X, \Delta)$ be a log Fano pair which admits a $\bT$-action. Then $(X, \Delta)$ is K-semistable if and only if $\Fut(\mX, \Delta_{\tc})\geq 0$ for any $\bT$-equivariant special test configuration $(\mX, \Delta_{\tc})\to \bA^1$ of $(X, \Delta)$.
\end{theorem}

\begin{proof}
The proof is the same as the second proof of Theorem \ref{thm: LX} (in Section \ref{subsec: fg}) up to replacing the divisorial minimizer $E$ with the $\bT$-equivariant divisorial minimizer given by Theorem \ref{thm: G-minimizer}.
\end{proof}

The following definition is inspired by Theorem \ref{thm: G-LX}.

\begin{definition}
Let $(X, \Delta)$ be a log Fano pair which admits a $\bT$-action. We say $(X, \Delta)$ is \textit{$\bT$-equivariant K-semistable} if $\Fut(\mX, \Delta_\tc)\geq 0$ for any $\bT$-equivariant special test configuration $(\mX, \Delta_\tc)$ of $(X, \Delta)$.
\end{definition}

\begin{corollary}\label{cor: kss=Tkss}
Let $(X, \Delta)$ be a log Fano pair which admits a $\bT$-action. Then $(X, \Delta)$ is K-semistable if and only if it is $\bT$-equivariant K-semistable.
\end{corollary}

\begin{proof}
 The result is implied by Theorem \ref{thm: G-LX}. Note that the result is well-known by \cite{LX20} when $(X, \Delta)$ is a log Fano $\bQ$-pair.
 \end{proof}

\subsection{Optimal degeneration}

In this subsection, we show that the divisor computing the $\delta$-invariant actually induces an optimal destabilization as in Definition \ref{def: optimal deg}.

\begin{theorem}\label{thm: optimal deg}
Let $(X, \Delta)$ be a log Fano pair with $\delta(X, \Delta)\leq 1$. Let $E$ be a prime divisor over $X$ computing $\delta(X, \Delta)$. Then $E$ induces a special test configuration of $(X, \Delta)$, denoted by $(\mX, \Delta_{\tc})\to \bA^1$, such that $\delta(X,\Delta)=\delta(\mX_0, \Delta_{\tc,0})$.
\end{theorem}

When $(X, \Delta)$ is a log Fano $\bQ$-pair, the above result is proved by \cite{BLZ22}. So we just assume $(X, \Delta)$ is a log Fano $\bR$-pair. We will take a similar strategy as in \cite{BLZ22}.

\begin{definition}\label{def: twisted kss}
Let $(X, \Delta)$ be a log Fano $\bR$-pair with $\delta(X, \Delta)\leq 1$. For a positive real number $0<\mu\leq 1$, we say $(X, \Delta)$ is \textit{$\mu$-twisted K-semistable} if $\Fut_{1-\mu}(\mX, \Delta_\tc)\geq 0$ for 
any special test configuration $(\mX, \Delta_{\tc})\to \bA^1$ of $(X, \Delta)$, where
$$\Fut_{1-\mu}(\mX, \Delta_\tc):=\sup_{D\in |-K_X-\Delta|_\bR}\Fut(\mX, \Delta_\tc+(1-\mu)\mD; \mL) \quad \text{and}\quad \mL:=-K_\mX-\Delta_\tc.$$
Note that $\mD$ is the natural extension of $D\in |-K_X-\Delta|_\bR$ and $(\mX, \Delta_\tc+(1-\mu)\mD; \mL)\to \bA^1$ is a test configuration of $(X, \Delta+(1-\mu)D)$.

Suppose $(X, \Delta)$ admits a $\bT:=\bG_m^r$-action for some $r\in \bN$. We say $(X, \Delta)$ is \textit{$\bT$-equivariant $\mu$-twisted K-semistable} if $\Fut_{1-\mu}(\mX, \Delta_\tc)\geq 0$ for 
any $\bT$-equivariant special test configuration $(\mX, \Delta_{\tc})\to \bA^1$ of $(X, \Delta)$.
\end{definition}

We list some properties under the above definition.

\begin{proposition}\label{prop: twisted kss}
Let $(X, \Delta)$ be a log Fano $\bR$-pair with $\delta(X, \Delta)\leq 1$. Fix a positive real number $0<\mu\leq 1$. Then we have the following conclusions:
\begin{enumerate}
\item for a non-trivial special test configuration $(\mX, \Delta_{\tc})\to \bA^1$ of $(X, \Delta)$, the sup in the definition of  $\Fut_{1-\mu}(\mX, \Delta_\tc)$ is achieved by a general element in $|-K_X-\Delta|_\bR$, moreover, 
$$\Fut_{1-\mu}(\mX, \Delta_\tc)=c\cdot (A_{X, \Delta}(E)-\mu\cdot S_{X, \Delta}(E)), $$
where $c\in \bQ_{>0}$ and $E$ is a prime divisor over $X$ with ${\ord_{\mX_0}}|_{K(X)}=c\cdot \ord_E$; 

\item $(X, \Delta)$ is $\mu$-twisted K-semistable if and only if $\delta(X, \Delta)\geq \mu$;

\item suppose $(X, \Delta)$ admits a $\bT:=\bG_m^r$-action for some $r\in \bN$, then $(X, \Delta)$ is $\mu$-twisted K-semistable if and only if
it is $\bT$-equivariant $\mu$-twisted K-semistable;

\item suppose $E$ is a prime divisor over $X$ computing $\delta(X, \Delta)$, then $E$ induces a special test configuration $(\mX, \Delta_\tc)\to \bA^1$ of $(X, \Delta)$ with $\Fut_{1-\delta(X, \Delta)}(\mX, \Delta_\tc)=0$; conversely, if $(\mX, \Delta_\tc)\to \bA^1$ is a special test configuration of $(X, \Delta)$ with $\Fut_{1-\delta(X, \Delta)}(\mX, \Delta_\tc)=0$, then the test configuration is induced by a prime divisor over $X$ computing $\delta(X, \Delta)$.
\end{enumerate}
\end{proposition}

\begin{proof}
The first statement is proved by the same way as \cite[Proposition 3.8]{BLZ22}, where the strategy does not depend on whether the coefficients are rational or not. Moreover, for any $D\in |-K_X-\Delta|_\bR$ whose support does not contain the center of $\ord_E$, we have
$$\Fut_{1-\mu}(\mX, \Delta_\tc)=\Fut(\mX, \Delta_\tc+(1-\mu)\mD; \mL)=c\cdot (A_{X, \Delta}(E)-\mu\cdot S_{X, \Delta}(E)).$$
Next we prove the statement (2). By the statement (1), $(X, \Delta)$ is $\mu$-twisted K-semistable if and only if $A_{X, \Delta}(E)-\mu\cdot S_{X, \Delta}(E)\geq 0$ for any special divisor $E$. By Theorem \ref{thm: special minimizer}, $\delta(X, \Delta)$ is computed by some special divisor, thus $\delta(X, \Delta)\geq \mu$ is equivalent to $A_{X, \Delta}(E)-\mu\cdot S_{X, \Delta}(E)\geq 0$ for any special divisor $E$.

We turn to the statement (3). Suppose $(X, \Delta)$ is $\bT$-equivariant $\mu$-twisted K-semistable.  If $\delta(X, \Delta)=1$, then $(X, \Delta)$ is automatically $\mu$-twisted K-semistable by statement (2). Assume $\delta(X, \Delta)<1$. By statement (1),  $A_{X, \Delta}(E)-\mu\cdot S_{X, \Delta}(E)\geq 0$ for any $\bT$-equivariant special divisor $E$. By Theorem \ref{thm: G-minimizer}, $\delta(X, \Delta)$ is computed by some $\bT$-equivariant special divisor, thus $\delta(X, \Delta)\geq \mu$ and $(X, \Delta)$ is $\mu$-twisted K-semistable by statement (2).

The last statement is a corollary of statement (1). 
\end{proof}

We are ready to prove Theorem \ref{thm: optimal deg}.

\begin{proof}[Proof of Theorem \ref{thm: optimal deg}]
By Theorem \ref{thm: special minimizer},  $E$ induces a special test configuration of $(X, \Delta)$, denoted by $(\mX, \Delta_{\tc})\to \bA^1$. Write $(X', \Delta'):=(\mX_0, \Delta_{\tc,0})$. It remains to show  $\delta(X,\Delta)=\delta(X',\Delta')$. 

We first show $\delta(X', \Delta')\geq \delta(X, \Delta)$. By Proposition \ref{prop: twisted kss} (2), it suffices to show that $(X', \Delta')$ is $\delta(X, \Delta)$-twisted K-semistable. Note that $(X', \Delta')$ admits a $\bT:=\bG_m$-action. Suppose $(X', \Delta')$ is not $\delta(X, \Delta)$-twisted K-semistable, by Proposition \ref{prop: twisted kss} (3), there exists a $\bT$-equivariant special test configuration of $(X', \Delta')$, denoted by $(\mX', \Delta'_{\tc})\to \bA^1$, such that 
$$\Fut_{1-\delta(X, \Delta)}(\mX', \Delta'_{\tc})<0.$$ 
On the other hand, by Proposition \ref{prop: twisted kss} (1), for general $D\in |-K_X-\Delta|_\bR$ we have
$$\Fut_{1-\delta(X, \Delta)}(\mX, \Delta_{\tc})=\Fut(\mX, \Delta_\tc+(1-\delta(X, \Delta))\mD; \mL)=0, $$
where $\mD$ is the natural extension of $D$ and $\mL:=-K_{\mX}-\Delta_\tc$.
By the proof of \cite[Lemma 3.1]{LWX21}, for $k\gg 1$, we can construct a new special test configuration $(\mX'', \Delta''_\tc)\to \bA^1$ of $(X, \Delta)$ such that for general $D\in |-K_X-\Delta|_\bR$ we have
\begin{align*}
&\ \Fut(\mX'', \Delta''_\tc+(1-\delta(X, \Delta))\mD'')\\
=&\ k\cdot \Fut(\mX, \Delta_\tc+(1-\delta(X, \Delta))\mD;\mL)+\Fut(\mX', \Delta'_\tc+(1-\delta(X, \Delta))\mD';\mL')\\
\leq&\  k\cdot \Fut_{1-\delta(X, \Delta)}(\mX, \Delta_{\tc})+\Fut_{1-\delta(X, \Delta)}(\mX', \Delta'_{\tc})\\
<&\ 0,
\end{align*}
where $\mD''$ is the natural extension of $D$ on $\mX''$, $\mD'$ is the natural extension of $\mD_0$ on $\mX'$,
and $\mL':=-K_{\mX'}-\Delta'_{\tc}$. This leads to a contradiction to the fact that $(X, \Delta)$ is $\delta(X, \Delta)$-twisted K-semistable. The contradiction implies that $(X', \Delta')$ is $\delta(X, \Delta)$-twisted K-semistable and $\delta(X', \Delta')\geq \delta(X, \Delta)$.

We next show $\delta(X', \Delta')\leq \delta(X, \Delta)$. Write 
$$\Delta:=\sum_{j=1}^kr_jD_j,\ \ \ \Delta':=\sum_{j=1}^kr_jD'_j, \ \quad\text{and}\quad \Delta_\tc:=\sum_{j=1}^kr_j\mD_j,$$ 
where $D_j$'s are effective Weil divisors. By Proposition \ref{prop: QR}, there exists a rational polytope $P$ containing $(r_1,...,r_k)$ as an interior point such that $(X, \sum_{j=1}^k x_j D_j)$ is a log Fano pair for any $(x_1,...,x_k)\in P$ and $(\mX, \sum_{j=1}^kx_j\mD_j)\to \bA^1$ is a special test configuration of $(X, \sum_{j=1}^kx_jD_j)$ for any $(x_1,...,x_k)\in P$.
For any rational vector $(x_1,...,x_k)\in P$, we have the following inequality by the lower semi-continuity of $\delta$-invariants (e.g. \cite{BL22}):
$$\delta(X', \sum_{j=1}^kx_jD_j')\leq \delta(X, \sum_{j=1}^kx_jD_j). $$
By the continuity of $\delta$-invariants (e.g. Proposition \ref{prop: continuity}), we obtain the desired $\delta(X', \Delta')\leq \delta(X, \Delta)$. The proof is complete.
\end{proof}

\section{Boundedness and Openness}\label{sec: stack}

Starting from this section, we will be devoted to construct a moduli space for K-semistable log Fano $\bR$-pairs by confirming several ingredients as in $\bQ$-coefficients case. In this section, we prove the boundedness and openness, and then construct a finite type Artin stack to parametrize K-semistable log Fano $\bR$-pairs with some fixed invariants.

\subsection{Boundedness}

Fixing a positive integer $d$, a positive real number $v$, and a finite set $I$ of positive real numbers. We construct a set $\mK:=\mK(d, v, I)$ of log Fano pairs as follows. We say a log Fano pair $(X, \Delta)$ is contained in $\mK$ if and only if the following conditions are satisfied:
\begin{enumerate}
\item $(X, \Delta)$ is a K-semistable log Fano pair of dimension $d$;
\item the coefficients of $\Delta$ are contained in $I$;
\item $\vol(-K_X-\Delta)\geq v$.
\end{enumerate}

\begin{theorem}\label{thm: bdd}
The set $\mK:=\mK(d, v, I)$ is log bounded.
\end{theorem}

Theorem \ref{thm: bdd} is well known for $\bQ$-coefficients case, and there are several proofs from \cite{Jiang20, Chen20, LLX19, XZ21} based on \cite{HMX14, Birkar19, Birkar21}. In the $\bR$-coefficients case, there is no difficulty to reduce the problem to the $\bQ$-coefficients case as in the proof of Corollary \ref{cor: bdd}. We provide a proof for the readers' convenience.

\begin{proof}[Proof of Theorem \ref{thm: bdd}]
Take $(X, \Delta)\in \mK$ and denote $\mu$ to be the minimal number in $I$. Suppose $\Delta:=\sum_{j=1}^kr_j D_j$ admits irrational coefficients, where $D_j$'s are effective Weil divisors. 
By Corollary \ref{cor: bdd}, $X$ lies in a bounded set. 
Thus there exists a very ample line bundle $A$ on $X$ with 
$$A^{d}\leq N_1\quad\text{and}\quad A^{d-1}\cdot (-K_X)\leq N_2,$$ 
where $N_1, N_2$ are two positive numbers which do not depend on the choice of $X$.
On the other hand, since $-K_X-\sum_{j=1}^kr_jD_j$ is pseudo-effective, we have
$$A^{d-1}\cdot D_j\leq \frac{A^{d-1}\cdot (-K_X)}{r_j}\leq \frac{N_2}{\mu},$$
which implies that $D_j$ is also bounded for each $j$. Thus $(X, \Delta)$ lies in a log bounded set. The proof is complete.
\end{proof}

\subsection{Openness}\label{subsec: openness}

Let $\pi: (Y, B)\to T$ be a family of $d$-dimensional log Fano pairs over a smooth base $T$ such that $Y\to T$ is flat. We say it is an \textit{$\bR$-Gorenstein family} if $-K_{Y/T}-B$ is $\bR$-Cartier. Fix a positive real number $1\leq a<\frac{d+1}{d}$. For each $t\in T$, we define $\tilde{\delta}_a(Y_t, B_t):=\min\{\delta(Y_t, B_t), a\}$.

\begin{theorem}\label{thm: constructible}
Let $\pi: (Y, B)\to T$ be an $\bR$-Gorenstein family of log Fano pairs over a smooth base $T$. Then the following function is constructible:
$$T\to \bR_{>0},\ \ t\mapsto \tilde{\delta}_a(Y_t, B_t). $$
\end{theorem}

\begin{proof}
By the proof of \cite[Theorem 5.4]{LZ24} (where we replace \cite[Proposition 5.1]{LZ24} with Proposition \ref{prop: qm}), there exists a constructible stratification $T=\amalg_i T_i$ such that for each index $i$, the invariant $\tilde{\delta}_a(Y_t, B_t)$ does not depend on $t\in T_i$.
\end{proof}

\begin{theorem}\label{thm: lsc}
Let $\pi: (Y, B)\to T$ be an $\bR$-Gorenstein family of log Fano pairs over a smooth base $T$. For any fixed $t_0\in T$, there exists an open neighborhood $t_0\in U\subset T$ such that $ \tilde{\delta}_a(Y_t, B_t)\geq \tilde{\delta}_a(Y_{t_0}, B_{t_0})$ for any $t\in U$.
\end{theorem}

\begin{proof}
Write $B=\sum_{j=1}^kr_jB_j$, where $B_j$'s are effective Weil divisors. If $(r_1,...,r_k)$ is rational, then the result is implied by \cite[Proposition 4.3]{BLX22} and \cite[Corollary 3.7]{LXZ22}. We assume $(r_1,...,r_k)$ is irrational. By Proposition \ref{prop: general QR}, there exists a rational polytope $P$ containing $(r_1,...,r_k)$ as an interior point such that $(Y, \sum_{j=1}^kx_jB_j)\to T$ is an $\bR$-Gorenstein family of log Fano pairs for any $(x_1,...,x_k)\in P$. Let $\{(r_{i1},...,r_{ik})\}_{i=1}^\infty\subset P$ be a sequence of rational vectors tending to $(r_1,...,r_k)$. For each index $i$, by \cite[Proposition 4.3]{BLX22} and \cite[Corollary 3.7]{LXZ22}, there exists an open neighborhood $t_0\in U_i$ such that 
$$\tilde{\delta}_a(Y_{t}, \sum_{j=1}^kr_{ij}B_{j, t})\geq \tilde{\delta}_a(Y_{t_0}, \sum_{j=1}^kr_{ij}B_{j, t_0})$$
for any $t\in U_i$. By the continuity, we have 
$$\tilde{\delta}_a(Y_t, B_t)\geq \tilde{\delta}_a(Y_{t_0}, B_{t_0})$$
for any $t\in \cap_{i=1}^\infty U_i$. The proof is finished by noting that  the function $ t\mapsto \tilde{\delta}_a(Y_t, B_t)$ is constructible (e.g.  Theorem \ref{thm: constructible}). 
\end{proof}

\begin{theorem}\label{thm: openness}
Let $\pi: (Y, B)\to T$ be an $\bR$-Gorenstein family of log Fano pairs over a smooth base $T$. Then the following sets are open subsets of $T$:
$$\{t\in T\ |\  \textit{$(Y_t, B_t)$ is K-semistable}\}\quad and \quad \{t\in T\ |\  \textit{$\delta(Y_t, B_t)>1$}\}.$$
\end{theorem}

\begin{proof}
Implied by Theorem \ref{thm: lsc}.
\end{proof}

\subsection{Moduli stack}\label{subsec: stack}

Given boundedness and openness, in this subsection, we construct a finite type Artin stack to parametrize K-semistable log Fano pairs in $\mK:=\mK(d, v, I)$.
Since $\mK$ is log bounded, there exist finite projective families of couples $(\mX^{(i)}, \mD^{(i)})\to T_i$ (indexed by $i$) over smooth basis $T_i$'s such that
\begin{enumerate}
\item for each $i$, $\mX^{(i)}\to T_i$ is proper and flat, and we could write $\mD^{(i)}=\sum_{j=1}^{k_i} a_{ij}\mD_{ij}$, where $\{a_{ij}\}_{j=1}^{k_i}\subset I$, $\mD_{ij}$'s are effective Weil divisors on $\mX^{(i)}$, and each $\mD_{ij}$ is flat over $T_i$ with no fiber contained in its support;
\item for any $(X, \Delta)\in \mK$, there is an index $i$ and a point $t\in T_i$ such that $(X, \Delta)$ is isomorphic to the fiber $(\mX^{(i)}_t, \mD^{(i)}_t)$.
\end{enumerate}
By \cite[Section 1.3]{LZ24}, we easily see the following set is finite:
$$\vol(\mK):=\{\vol(-K_X-\Delta)\ |\ \textit{$(X, \Delta)\in \mK$}\}.$$
Fixing a real number $v_0\in \vol(\mK)$ and an \textit{irrational} vector $\fa:=(a_1,...,a_k)$, where $\{a_j\}_{j=1}^k\subset I$, we aim to construct a finite type Artin stack parametrizing the following set:
$$\mK_{d, v_0, \fa}:=\{(X, \sum_{j=1}^k a_jD_j)\in \mK\ |\ \textit{$D_j$'s are Weil divisors and $\vol(-K_X-\sum_{j=1}^k a_jD_j)=v_0$}\}. $$

\begin{proposition}\label{prop: epsilon-delta}
There exist a rational polytope $P\subset \bR^k_{>0}$ containing $\fa$ as an interior point and a positive real number $0<\epsilon_0<1$ such that the following condition is satisfied: for any $(X, \sum_{j=1}^k a_jD_j)\in \mK_{d,v_0,\fa}$, $(X, \sum_{j=1}^k x_jD_j)$ is a log Fano pair with $\vol(-K_X-\sum_{j=1}^k x_jD_j)\geq \epsilon_0$ for any $(x_1,...,x_k)\in P$.
\end{proposition}

\begin{proof}
Note that $\mK_{d, v_0, \fa}$ is a log bounded set of K-semistable log Fano pairs with the fixed volume $v_0$. By \cite[Propositions 4.1, 3.5, 3.3]{LZ24}, there exists a rational polytope $P$ containing $\fa$ as an interior point such that $(X, \sum_{j=1}^kx_jD_j)$ is log Fano for any $(X, \sum_{j=1}^k a_jD_j)\in \mK_{d,v_0,\fa}$ and any $(x_1,...,x_k)\in P$.
By the continuity of the volume, we easily get the desired $\epsilon_0$ up to shrinking $P$.
\end{proof}

Note that $P$ does not depend on the choice of $(X, \sum_{j=1}^k a_jD_j)\in \mK_{d,v_0,\fa}$. By Koll\'ar's notation, $P$ is a rational polytope in the linear $\bQ$-envelope $\LEnv_\bQ(1, \fa)$ (e.g. \cite[Definition 8.11, 11.44]{Kollar23}). 
We are ready to define the moduli functor $\mM^{\Kss}_{d, v_0, \fa}$. Note that $\fa$ is irrational.
For any scheme $S$, a family $(\mX, \sum_{j=1}^ka_j\mD_j)\to S$ is contained in $\mM^\Kss_{d,v_0,\fa}(S)$ if and only if the following conditions are satisfied:
\begin{enumerate}
\item $\mX\to S$ is proper and flat;
\item each $\mD_j$ is a K-flat family of Weil divisors on $\mX$ (e.g. \cite[Chapter 7]{Kollar23});
\item $-K_{\mX/S}-\sum_{j=1}^k x_j\mD_j$ is $\bR$-Cartier for any $(x_1,...,x_k)\in P$;
\item for each closed point $s\in S$, the fiber $(\mX_s, \sum_{j=1}^ka_j\mD_{j,s})$ is contained in $\mK_{d, v_0, \fa}$.
\end{enumerate}
Note that in the third condition, $-K_{\mX/S}-\sum_{j=1}^k x_j\mD_j$ is automatically relatively ample by the construction of $P$.  The third condition is a little different from that in \cite[Section 2.6]{XZ20b} since we are dealing with irrational coefficients. The definition here is in the same spirit of \cite[8.13]{Kollar23} (see \cite[4.20]{Kollar23}). However, if we restrict to reduced bases, the third condition above could be replaced with the following standard one due to Proposition \ref{prop: U} and \cite[Corollary 4.35]{Kollar23}:

\ (3') $-K_{\mX/S}-\sum_{j=1}^k a_j\mD_j$ is $\bR$-Cartier.

\begin{proposition}\label{prop: U}
Notation as in Proposition \ref{prop: epsilon-delta}. There exists a finite rational polytope chamber decomposition $P:=\cup_s P_s$ such that for any face $F$ of any chamber $P_s$ and for any family of couples $(\mX, \sum_{j=1}^k\mD_j)\to S$ satisfying
\begin{enumerate}
\item $\mX\to S$ is proper and flat over a smooth base $S$,
\item each $\mD_j$ is an effective Weil divisor on $\mX$ with no fiber contained in its support, 
\item $(\mX_s, \sum_{j=1}^ka_j\mD_{j,s})$ is contained in $\mK_{d,v_0,\fa}$ for any $s\in S$,
\end{enumerate}
we have the following conclusion: if 
$-K_{\mX/S}-\sum_{j=1}^kc_j\mD_j$ is $\bR$-Cartier for some $(c_1,...,c_k)\in F^\circ$, then $-K_{\mX/S}-\sum_{j=1}^kx_j\mD_j$ is $\bR$-Cartier for any $(x_1,...,x_k)\in F^\circ$.
\end{proposition}

\begin{proof}
Define the set $\mP:=\{(X, \sum_{j=1}^k D_j)\ |\ \textit{$(X, \sum_{j=1}^ka_jD_j)\in \mK_{d, v_0, \fa}$}\}$. By Proposition \ref{prop: epsilon-delta}, the set $\mP$ is contained in $\mA_{1}(1, \epsilon_0)$ under the notation in \cite[Section 7]{LZ24}.
The result is thus implied by \cite[Proposition 7.4]{LZ24}. 
\end{proof}

It is clear that there must be a chamber $P_s$ containing $\fa$ as an interior point. Thus if we take $F=P_s$, the conclusion in Proposition \ref{prop: U} holds for $P_s$.

We aim to show that the moduli functor $\mM^\Kss_{d, v_0, \fa}$ is a finite type Artin stack.
\begin{theorem}\label{thm: artin stack}
The moduli functor $\mM^\Kss_{d,v_0,\fa}$ is a finite type Artin stack.
\end{theorem}

\begin{proof}
Since $\mK_{d, v_0, \fa}$ is log bounded, there exist finite projective families of couples $(\mX^{(i)}, \mD^{(i)})\to T_i$ (indexed by $i$) over smooth basis $T_i$'s such that
\begin{enumerate}
\item for each $i$, $\mX^{(i)}\to T_i$ is proper and flat, and we could write $\mD^{(i)}:=\sum_{j=1}^{k} \mD_{ij}$, where $\mD_{ij}$'s are effective Weil divisors on $\mX^{(i)}$ and each $\mD_{ij}$ is flat over $T_i$ with no fiber contained in its support;
\item for any $(X, \sum_{j=1}^k a_j D_j)\in \mK_{d, v_0, \fa}$, there is an index $i$ and a point $t\in T_i$ such that $(X, \sum_{j=1}^k D_j)$ is isomorphic to the fiber $(\mX^{(i)}_t, \mD^{(i)}_t)$.
\end{enumerate}
By \cite[Section 1.3]{LZ24}, up to a further constructible stratification of each $T_i$, we may assume the log Fano domain $\LF(\mX_t^{(i)}, \mD^{(i)}_t)$ (see \cite[Definition 3.1]{LZ24}) does not depend on $t\in T_i$ for each $i$. By Proposition \ref{prop: epsilon-delta}, we know that $P\subset \LF(X, \sum_{j=1}^k D_j)$ for any $(X, \sum_{j=1}^k a_j D_j)\in \mK_{d, v_0, \fa}$. Up to subtracting redundant families, we may assume that $P\subset \LF(\mX_t^{(i)}, \mD^{(i)}_t)$ for any index $i$ and any $t\in T_i$. By the same argument as in the first paragraph of \cite[Proof of Proposition 3.5]{LZ24}, up to a further constructible stratification of each $T_i$, we may assume $-K_{\mX^{(i)}/T_i}-\sum_{j=1}^k x_j \mD_{ij}$ is $\bR$-Cartier for any index $i$ and any $(x_1,...,x_k)\in P$. This implies that the Hilbert function 
$$m\  \text{(sufficiently divisible)}\mapsto \chi_i(t; (x_1,...,x_k)):=\chi\left(-m(K_{\mX_t^{(i)}}+\sum_{j=1}^k x_j\mD_{ij, t})\right)$$
does not depend on $t\in T_i$ for any fixed rational vector $(x_1,...,x_k)\in P$. Let $P(\bQ)$ be the set consisting of rational vectors in $P$. For each $(X, \sum_{j=1}^k a_j D_j)\in \mK_{d, v_0, \fa}$, we could associate a function defined on $P(\bQ)\times \bZ_+$\footnote{For any rational vector $(x_1,...,x_k)\in P$, we only consider sufficiently divisible $m$ for the second factor.} as follows:
$$(x_1,...,x_k)\times m\mapsto \chi\left(-m(K_X+\sum_{j=1}^kx_jD_j)\right). $$
Thus there are only finitely many such functions, denoted by
$$h_l((x_1,...,x_k), m), \ \ \ 1\leq l\leq s. $$
Let $\mM_l\subset \mM^\Kss_{d, v_0, \fa}$ be the sub-functor parametrizing objects whose associated function is $h_l$. Thus we see $\mM^\Kss_{d, v_0,\fa}=\coprod_l \mM_l$. It suffices to show that $\mM_l$ is a finite type Artin stack. 

Let $\{\fa_r:=(a_{r1},...,a_{rk})\}_r$ be the finite set of vertices of $P$. For each $r$, we choose a sufficiently divisible positive integer $N_r$ such that $-N_r(K_X+\sum_{j=1}^k a_{rj}D_j)$ is a very ample line bundle for any $(X, \sum_{j=1}^k a_jD_j)\in \mM_l(\bC)$. Consider the following set of vectors of degrees:
$$S_r:=\{(-N_r(K_X+\sum_{j=1}^k a_{rj}D_j)\cdot D_1,...,-N_r(K_X+\sum_{j=1}^k a_{rj}D_j)\cdot D_k) | \textit{$(X, \sum_{j=1}^k a_jD_j)\in \mM_l(\bC)$}\}. $$
It is clear that for each $r$, the set $S_r$ is finite. Set $M:=h_l((a_{11},...,a_{1k}), N_1)-1$. Let $\Hilb_l(\bP^M)$ be the Hilbert scheme parametrizing closed subschemes of $\bP^M$ with Hilbert polynomial $h_l((a_{11},...,a_{1k}), mN_1)$. Let $U\subset \Hilb_l(\bP^M)$ be the open subscheme parametrizing normal and Cohen-Macaulay varieties. By \cite[Theorem 7.3]{Kollar23}, there exists a separated $U$-scheme $W_1$ of finite type parametrizing K-flat divisors $(D_1,...,D_k)$ with all possible vectors of degrees in $S_1$. Write $(\fX, \sum_{j=1}^k \fD_j)\to W_1$ for the corresponding universal family. It is clear that $\mM_l(\bC)$ is contained in the set of fibers of $(\fX, \sum_{j=1}^k a_j\fD_j)\to W_1$. Denote $L_r:=\omega_{\fX/W_1}^{[-N_r]}(-N_r\cdot \sum_{j=1}^k a_{rj}\fD_j)$. By \cite[Corollary 3.30]{Kollar23}, there exists a locally closed subscheme $W_2\subset W_1$ such that a map $T\to W_1$ factors through $W_2$ if and only if 
$$(L_r)_T:= \omega_{\fX_T/T}^{[-N_r]}(-N_r\cdot \sum_{j=1}^k a_{rj}(\fD_j)_T)$$
is invertible for any index $r$, where $(\fX_T, \sum_{j=1}^k (\fD_j)_T)\to T$
is obtained by the base change along $T\to W_1$. By Theorem \ref{thm: openness}, the following set is open in $W_2$:
$$W:=\{t\in W_2\ |\  \textit{$(\fX_t, \sum_{j=1}^k a_j \fD_{j,t})$ is a K-semistable log Fano pair}\}. $$
By the above construction, we see that $\mM_l\cong [W/\PGL(M+1)]$, which is indeed a finite type Artin stack. The proof is complete.
\end{proof}

\section{$\Theta$-reductivity and K-polystable degeneration}\label{sec: Theta}

In this section, we aim to confirm the $\Theta$-reductivity of $\mM^\Kss_{d, v_0, \fa}$ (see Section \ref{subsec: stack}), and show that there exists a unique K-polystable degeneration of a given K-semistable log Fano $\bR$-pair. Throughout this section, we fix $R$ to be a DVR essentially of finite type over $\bC$ with fraction field $K$ and residue field $\kappa$. We write $t$ to be the parameter of the affine line $\bA^1$. Let $0_K\in \bA^1_K$ be the closed point defined by $t=0$, and $0_{\kappa}\in \bA^1_R$ the one defined by the vanishing of $t$ and a uniformizing parameter $\pi\in R$. Let $0_R$ be the line in $\bA^1_R$ defined by $t=0$.

\subsection{$\Theta$-reductivity}

\begin{definition}
Denote $\Theta:=[\bA^1/\bG_m]$. We write $\Theta_R:=\Theta\times \Spec R\cong [\Spec R[t]/\bG_m]$, which is a stacky surface with the unique closed point $0_\kappa$. We say $\mM^\Kss_{d, v_0, \fa}$ is $\Theta$-reductive if for any DVR (essentially of finite type) $R$ and any map $\Theta_R\setminus 0_\kappa\to \mM^\Kss_{d, v_0, \fa}$, there exists a unique extension $\Theta_R\to \mM^\Kss_{d, v_0, \fa}$. 
\end{definition}

Denote $\Theta_K:=[\Spec K[t]/\bG_m]$, which is an open sub-stack of $\Theta_R$ defined by $\pi\ne 0$. It is clear that the map $\phi^\circ:\Theta_R\setminus 0_\kappa\to \mM^\Kss_{d, v_0, \fa}$ is equivalent to a map $\Spec R\to \mM^\Kss_{d, v_0, \fa}$ together with a map $\phi': \Theta_K\to \mM^\Kss_{d, v_0, \fa}$ with $\phi^\circ|_{\Theta_K}=\phi'$.
We will only assume $\fa=(a_1,...,a_k)$ is irrational.

\begin{theorem}\label{thm: Theta}
Given $[(X, D(\fa):=\sum_{j=1}^ka_jD_j)\to  \Spec R]\in \mM^\Kss_{d, v_0, \fa}$ and a map $\phi^\circ: \Theta_R\setminus 0_\kappa\to \mM^\Kss_{d, v_0, \fa}$. Then there exists a unique extension $\phi: \Theta_R\to \mM^\Kss_{d, v_0, \fa}$. In particular, $\mM^\Kss_{d, v_0, \fa}$ is $\Theta$-reductive.
\end{theorem}

\begin{proof}
We divide the proof into several steps.

\

\textbf{Step 1}. \textit{In this step, we show that $\phi^\circ|_{\Theta_K}$ induces a special test configuration of $(X_K, D(\fa)_K)$ with vanishing generalized Futaki invariant.} 
It is clear that $\phi^\circ|_{\Theta_K}$ induces a special test configuration of $(X_K, D(\fa)_K)$, denoted by
$$(\mX_K, \mD(\fa)_K:=\sum_{j=1}^k a_j \mD_{j, K})\to \bA^1_K.$$ 
Since the central fiber $(\mX_{0_K}, \mD(\fa)_{0_K})$ is K-semistable, 
we see $\Fut(\mX_K, \mD(\fa)_K)=0$ by Corollary \ref{cor: kss deg}. 

If $\delta(\mX_K, \mD(\fa)_K)>1$, then $(\mX_K, \mD(\fa)_K)$ is a trivial test configuration, otherwise there is a contradiction by applying Proposition \ref{prop: fut=beta}. In this case $\phi^\circ$ naturally extends cross $0_\kappa$. From now on, we assume $\delta(\mX_K, \mD(\fa)_K)=1$ and $(\mX_K, \mD(\fa)_K)$ is a non-trivial special test configuration. We write $(\mX^\circ, \mD^\circ)\to \bA^1_R\setminus 0_\kappa$ for the corresponding family with respect to $\phi^\circ$.

\

\textbf{Step 2}. \textit{In this step, we create a prime divisor $F$ over $X$.} By step $1$, there exists a special divisorial valuation $c\cdot \ord_{F_K}$ corresponding to the special test configuration $(\mX_K, \mD(\fa)_K)$, where $F_K$ is a prime divisor over $X_K$. Thus there exists a prime divisor $F$ over $X$ whose restriction over $X_K$ is exactly $F_K$. We aim to find an effective $\bR$-divisor $\Gamma\sim_\bR  -K_X-D(\fa)$ on $X$ such that $(X, D(\fa)+\Gamma)$ is lc with $F$ being an lc place. 

By Proposition \ref{prop: twisted kss} (4), $F_K$ computes $\delta(X_K, D(\fa)_K)=1$, i.e. $A_{X_K, D(\fa)_K}(F_K)-S_{X_K, D(\fa)_K}(F_K)=0$. Take a sequence of positive numbers $\{\delta_m\}_m$ satisfying 
$$\delta_m<\min\{\delta_m(X_K, D(\fa)_K), \delta_m(X_\kappa, D(\fa)_\kappa), 1\}\quad \text{and}\quad   \lim_m \delta_m=1.$$
For a large positive integer $m$, we take an $m$-basis type divisor $B_{m,K}$ for $(X_K, D(\fa)_K)$ such that $\ord_{F_K}(B_{m,K})=S_m(F_K)$. Observing the following
$$A_{X_K, D(\fa)_K}(F_K)-\delta_m\cdot S_m(F_K) \to 0 \ \ \text{as}\  \ m\to \infty,$$
we see that for any $0<\epsilon<1$, the following holds for $m\gg 1$:
$$0<A_{X_K, D(\fa)_K+\delta_m B_{m, K}}(F_K)<\epsilon.$$
Let $B_m$ be the closure of $B_{m, K}$ on $X$. 
Recalling that $(X_\kappa, D(\fa)_\kappa)$ is also K-semistable and $B_{m, K}$ degenerates to an $m$-basis type divisor $B_{m, \kappa}$ for $(X_\kappa, D(\fa)_\kappa)$, we see that $(X, D(\fa)+X_\kappa+\delta_m B_m)$ is lc with 
$$0<A_{X, D(\fa)+X_\kappa+\delta_m B_{m}}(F)<\epsilon. $$
Thus for any $0<\epsilon\ll 1$, one can easily cook up an effective $\bR$-divisor $\Gamma\in |-K_X-D(\fa)|_\bR$ such that $(X, D(\fa)+X_\kappa+\Gamma)$ is lc with $A_{X, D(\fa)+X_\kappa+\Gamma}(F)<\epsilon$. By a similar argument as in the proof of \cite[Lemma 3.2]{LXZ22}, there exists an effective $\bR$-divisor $\Gamma'\in |-K_X-D(\fa)|_\bR$ such that $(X, D(\fa)+X_\kappa+\Gamma')$ is lc with $F$ being its lc place.

\

\textbf{Step 3}. \textit{In this step, we construct the extending family.} 
By step 2, there exists an effective $\bR$-divisor $\Gamma\sim_\bR -K_X-D(\fa)$ on $X$ such that $(X, D(\fa)+X_\kappa+\Gamma)$ is lc with $F$ being its lc place. We aim to use $F$ to construct a degeneration of $(X, D(\fa)+X_\kappa+\Gamma)$. By \cite{BCHM10}, there exists a birational model (over $\Spec R$) which extracts precisely $F$, denoted by
$$ (Y, \widetilde{D(\fa)}+Y_\kappa+\widetilde{\Gamma}+F)\to (X, D(\fa)+X_\kappa+\Gamma), $$
where $\widetilde{D(\fa)}+\widetilde{\Gamma}$ is the birational transform of $D(\fa)+ \Gamma$. Consider the family 
$$(Y\times \bA^1, (\widetilde{D(\fa)}+Y_\kappa +\widetilde{\Gamma})\times \bA^1+F\times \bA^1+Y\times 0)\to \bA^1_R.$$
Denote $\Lambda:=\widetilde{D(\fa)}+Y_\kappa +\widetilde{\Gamma}$. Note that $F\times \bA^1$ and $Y\times 0$ are both lc places of 
$$(Y\times \bA^1, \Lambda\times \bA^1+F\times \bA^1+Y\times 0).$$
Let $b_1, b_2$ be two positive integers. Considering the weighted blowup with weight $(b_1,b_2)$ respectively along $F\times \bA^1$ and $Y\times 0$, we get a crepant  morphism
$$(\mY, \Lambda_\mY+F_\mY+\mY_{0_R}+\mE)\to  (Y\times \bA^1, \Lambda\times \bA^1+F\times \bA^1+Y\times 0),$$
where $\Lambda_\mY+ F_\mY+ \mY_{0_R}$ is the birational transform of $\Lambda\times \bA^1+F\times \bA^1+ Y\times 0$, and $\mE$ is the exceptional divisor centered in $Y\times 0$. Note that $\mY$ is Fano type over $\bA^1_R$ and $(\mY, \Lambda_\mY+F_\mY+\mY_{0_R}+\mE)$ admits a $\bG_m$-action which lifts the $\bG_m$- action on $\bA^1$. Up to a $\bG_m$-equivariant small $\bQ$-factorization, we may assume $\mY$ is $\bQ$-factorial. 
Thus one could run $\bG_m$-equivariant MMP/$\bA^1_R$ on $F_\mY+\mY_{0_R}$ which contracts $F_\mY+\mY_{0_R}$, and we denote the resulting model to be
$$(\mY, \Lambda_\mY+F_\mY+\mY_{0_R}+\mE)\dashrightarrow (\mY', \Lambda_{\mY'}+\mE'), $$
where $\Lambda_{\mY'}+\mE'$ are the push-forward of $\Lambda_\mY+\mE$. Writing $\mD_{\mY'}+\mY'_\kappa+\Gamma_{\mY'}$ to be the birational transform of $\widetilde{D(\fa)}\times \bA^1+Y_\kappa\times \bA^1+\widetilde{\Gamma}\times \bA^1$ on $\mY'$, we see
$$ \Lambda_{\mY'}=\mD_{\mY'}+\mY'_\kappa+\Gamma_{\mY'}.$$
Running $\bG_m$-equivariant MMP on $-K_{\mY'}-\mD_{\mY'}$ to get its anti-canonical model, denoted by 
$$(\mY', \mD_{\mY'}+\mY'_\kappa+\Gamma_{\mY'}+\mE')\dashrightarrow (\mY'', \mD_{\mY''}+\mY''_\kappa+\Gamma_{\mY''}+\mE''),$$
where $\mD_{\mY''}+\mY''_\kappa+\Gamma_{\mY''}+\mE''$ is the push-forward of $\mD_{\mY'}+\mY'_\kappa+\Gamma_{\mY'}+\mE'$.  Now we get a morphism 
$$(\mY'', \mD_{\mY''}+\mY''_\kappa+\Gamma_{\mY''}+\mE'')\to \bA^1_R$$
with $\mY''_{0_R}=\mE''$. 
By our construction, $(\mY'', \mD_{\mY''})\times \bA^1_K$ is a test configuration of $(X_K, D(\fa)_K)$ with the integral central fiber $\mY''_{0_K}=\mE''_{0_K}$. By the proof of \cite[Theorem 4.6]{BHJ17}, we can choose $(b_1, b_2)$ such that the restriction of $\ord_{\mE''_{0_K}}$ to the function field of $X_K$ is exactly $c\cdot \ord_{F_K}$ (see step 1). Since $-K_{\mY''}-\mD_{\mY''}$ is ample over $\bA^1_R$, we see that $(\mY'', \mD_{\mY''})$ is indeed isomorphic to $(\mX_K, D(\fa)_K)$ over $\bA^1_K$, and as a consequence $(\mY'', \mD_{\mY''})\to \bA^1_R$ is an extension of $\phi^\circ: (\mX^\circ, \mD^\circ)\to \bA^1_R\setminus 0_\kappa$. 
On the other hand, $(\mY'', \mD_{\mY''}+\mY''_\kappa+\mY''_{0_R}+\Gamma_{\mY''})$ is lc and 
$$K_{\mY''}+\mD_{\mY''}+\mY''_\kappa+\mY''_{0_R}+\Gamma_{\mY''}\sim_{\bR,  \bA^1_R} 0.$$
By adjunction, we see that 
$$(\mY''_\kappa, \mD_{\mY''_\kappa}+\mY''_{0_\kappa}+\Gamma_{\mY''_\kappa}):=(\mY'', \mD_{\mY''}+\mY''_\kappa+\mY''_{0_R}+\Gamma_{\mY''})|_{\mY''_\kappa}$$
is log canonical, and $(\mY''_\kappa, \mD_{\mY''_\kappa})\to \bA^1_\kappa$ is a weakly special test configuration of $(X_\kappa, D(\fa)_\kappa)$.

\

\textbf{Step 4}. \textit{In this step, we finish the proof}. By step 3, we obtain an extension of $\phi^\circ: (\mX^\circ, \mD^\circ)\to \bA^1_R\setminus 0_\kappa$, denoted by $\phi: (\mX, \mD)\to \bA^1_R$, and $(\mX_\kappa, \mD_\kappa)\to \bA^1_\kappa$ is a weakly special test configuration of $(X_\kappa, D(\fa)_\kappa)$. Observing the following relation (see step 1):
$$\Fut(\mX_K, \mD_K)=\Fut(\mX_\kappa, \mD_\kappa)=0, $$
we see that $(\mX_\kappa, \mD_\kappa)\to \bA^1_\kappa$ is actually a special test configuration of $(X_\kappa, \mD_\kappa)$ by Proposition \ref{prop: weak LX}. Applying Proposition \ref{prop: twisted kss} (4) and Theorem \ref{thm: optimal deg}, we see that the central fiber $(\mX_{0_\kappa}, \mD_{0_\kappa})$ is a K-semistable log Fano pair. This implies that $\Theta_R\setminus 0_\kappa\to \mM^\Kss_{d, v_0, \fa}$ (which corresponds to $\phi^\circ$) indeed extends to $\Theta_R\to \mM^\Kss_{d, v_0, \fa}$ (which corresponds to $\phi$).

Finally we show this extension is unique. Write the extension in the following way: $(\mX, \mD(\fa):=\sum_{j=1}^ka_j\mD_j)\to \bA^1_R$. By our construction, it is a flat family of log Fano $\bR$-pairs with $-K_{\mX/\bA^1_R}-\mD(\fa)$ being $\bR$-Cartier. By Proposition \ref{prop: general QR}, we could replace $\fa$ with a rational vector $\fa':=(a_1',...,a_r')$ such that $(\mX, \mD(\fa'))\to \bA^1_R$ is a flat family of log Fano $\bQ$-pairs with $-K_{\mX/\bA^1_R}-\mD(\fa')$ being $\bQ$-Cartier. Denote $\mL^\circ:=-N(K_{\mX/\bA^1_R}+\mD(\fa'))|_{\mX^\circ}$ and $\mL:=-N(K_{\mX/\bA^1_R}+\mD(\fa'))$, where $N$ is a sufficiently divisible positive integer such that $\mL^\circ$ is a line bundle. Applying \cite[Lemma 2.16]{ABHLX20}, we see the extension is indeed unique. The proof is complete.
\end{proof}

\subsection{K-polystable degeneration}

In this subsection, we aim to prove the following result on the K-polystable degeneration for a K-semistable log Fano $\bR$-pair.

\begin{theorem}\label{thm: kps deg}
Let $(X, \Delta)$ be a K-semistable log Fano $\bR$-pair. Then it admits a unique K-polystable degeneration via a special test configuration with vanishing generalized Futaki invariant. 
\end{theorem}

For a log Fano $\bQ$-pair, the theorem is well-known by \cite{LWX21}. Armed with the $\Theta$-reductivity in Theorem \ref{thm: Theta}, the proof is standard (e.g. \cite[Section 3.2]{LWX21} and \cite[Theorem 5.2]{BLZ22}). We sketch it for the readers' convenience.

\begin{proof}[Sketch proof of Theorem \ref{thm: kps deg}]
Suppose $\delta(X, \Delta)>1$, then any special test configuration of $(X, \Delta)$ with a K-semistable central fiber must be trivial by the same argument as in step 1 of the proof of Theorem \ref{thm: Theta}. This means $(X, \Delta)$ is K-polystable and any K-polystable degeneration is itself achieved by a trivial test configuration.

We assume $\delta(X, \Delta)=1$ and let $(\mX^{(i)}, \Delta^{(i)}_\tc)\to \bA^1$ ($i=1,2$) be two special test configurations with K-semistable central fibers. Applying Theorem \ref{thm: Theta} and Theorem \ref{thm: optimal deg}, one can easily cook up a family $(\fX, \fD)\to \bA_{t_1}^1\times \bA_{t_2}^1$ such that $(\fX_{0,0}, \fD_{0,0})$ is a common degeneration of $(\mX^{(i)}_0, \Delta^{(i)}_{\tc,0})$ for $i=1,2$ by two special test configurations, and $(\fX_{0,0}, \fD_{0,0})$ is K-semistable (see \cite[Theorem 5.5]{BLZ22}). Arguing by the same way as 
\cite[Section 3.2, Proof of theorem 1.3]{LWX21}, we know that there must be a K-polystable degeneration. Arguing by the same way as \cite[Proposition 5.7]{BLZ22}, the K-polystable degeneration is unique (up to isomorphism). The generalized Futaki invariant vanishes by Corollary \ref{cor: kss deg}. 
\end{proof}

\section{$S$-completeness}\label{sec: S}

In this section, we aim to show that the Artin stack $\mM^{\Kss}_{d, v_0, \fa}$ (see Section \ref{subsec: stack}) is $S$-complete (see Section \ref{subsec: S}). The essential idea is the same as \cite[Section 5]{BX19}, where the language of filtrations is needed. However, the different point here is that it seems subtle to directly define a filtration on the anti-log canonical ring of a given log Fano $\bR$-pair and show the finite generation of the graded algebra associated to the filtration. Our strategy here is still to conduct an approximation process. This means we only need the language of filtrations for log Fano $\bQ$-pairs.

\subsection{Filtrations}
In this subsection, we fix $(X, \Delta)$ to be a log Fano $\bQ$-pair and write $L:=-K_X-\Delta$. Denote by $R:=\oplus_{m\in M(L)} H^0(X, mL)$, where $M(L)$ consists of natural numbers divisible by the Cartier index of $L$.

\begin{definition}\label{def: filtration}
By a \emph{filtration} $\mF$ of $R$, we mean the data of a family of $\mathbb{C}$-vector subspaces
$\mF^\lambda R_m\subset R_m$
for $m\in M(L)$ and $\lambda\in \bR$, satisfying:
\begin{enumerate}
\item (decreasing) $\mF^\lambda R_m\subset \mF^{\lambda'}R_m$ if $\lambda\geq \lambda'$;
\item (left continuous) $\mF^\lambda R_m=\cap_{\lambda'<\lambda}\mF^{\lambda'}R_m$ for all $\lambda\in \bR$;
\item (linearly bounded) there exist $e_-, e_+\in \bR$ such that $\mF^{mx}R_m=0$ for all $x\geq e_+$ and $\mF^{mx}R_m=R_m$ for all $x\leq e_-$;
\item (multiplicative) $\mF^\lambda R_m\cdot \mF^{\lambda'}R_{m'}\subset \mF^{\lambda+\lambda'}R_{m+m'}$.
\end{enumerate}
\end{definition}

Given a filtration $\mF$ on $R$. For $m\in M(L)$ and $\lambda\in \bR$ we set
$$I_{m,\lambda}= I_{m,\lambda}(\mF):={\rm Im}\left(\mF^\lambda R_m\otimes \mO_X\left(-mL\right)\to \mO_X \right),$$
where the map is naturally induced by the evaluation $H^0(X,mL)\otimes \mO_X(-mL)\to \mO_X$. Define $I_m^{(t)}=I_m^{(t)}(\mF):=I_{m,mt}(\mF)$ for $t\in \bR$, then $I_\bullet^{(t)}$ is a graded sequence of ideals on $X$. For each fixed $\lambda\geq 0$, by Definition \ref{def: filtration} (4), we see that $\{I_{m, \lambda}\}_{m\in M(L)}$ is an increasing sequence of ideals. We write $\kb_\lambda(\mF):=\max_{m\in M(L)} I_{m, \lambda}$ and call $\kb_\bullet(\mF):=\{\kb_p(\mF)\}_{p\in \bN}$ \textit{the base ideal sequence} associated to $\mF$.

\begin{definition}{\rm (\cite[Definition 4.1]{XZ20b})}\label{def: slope}
The \emph{log canonical slope} associated to $\mF$ is defined as follows:
$$\mu_{X, \Delta}(\mF):=\sup \left\{t\in \bR\mid \textit{$\lct\left(X,\Delta; I_\bullet^{(t)}\right)\geq 1$}\right\}. $$
Following \cite{XZ20b}, we define the following beta-invariant associated to $\mF$:
$$\beta_{X, \Delta}(\mF):=\mu_{X, \Delta}(\mF)-S_{X, \Delta}(\mF), $$
where $S_{X, \Delta}(\mF)$ is the standard $S$-invariant associated to the filtration $\mF$ (see \cite{BJ20, AZ22}).
\end{definition}

\begin{definition}\label{def: ding}{\rm (\cite[Definition 2.13]{XZ20b})}
Let $\mF$ be a filtration on $R$ and choose $e_-,e_+\in \bZ$ as in Definition \ref{def: filtration}. Let $e:=e_+-e_-$. For each $m\in M(L)$, we set
$$\mI_m=\mI_m(\mF):= I_{m,me_+}+I_{m, me_+-1}\cdot t+...+I_{m, me_-+1}\cdot t^{me-1}+t^{me}\subset \mO_{X\times \bA^1}.$$
It is not hard to check that $\mI_\bullet$ is a graded sequence of ideals. Write $X_0:=X\times 0\subset X\times \bA^1$. We set
$$c_m:=\lct(X\times \bA^1, (\Delta\times \bA^1)\cdot (\mI_m)^{\frac{1}{m}}; X_0) \quad \text{and} \quad c_\infty:=\lim_{m\to \infty} c_m.$$
We define two invariants associated to $\mF$ as follows:
$$\fL^{\NA}(\mF):=c_\infty+e_+-1,\ \ \ \fD^{\NA}(\mF):=\fL^{\NA}(\mF)-S_{X, \Delta}(\mF).$$
\end{definition}

The following lemmas will be applied in the next subsection.

\begin{lemma}\label{lem: slope-lct}
Let $\mF$ be a filtration on $R$. Then $\mu_{X, \Delta}(\mF)\leq \lct(X, \Delta; \kb_\bullet(\mF))$.
\end{lemma}

\begin{proof}
For any positive rational number $a<\mu(\mF)$, we see that $\lct(X, \Delta; I_\bullet^{(a)}(\mF))\geq 1$. For any sufficiently divisible $m$, we have $I_{m, ma}(\mF)\subset \kb_{ma}(\mF)$, thus 
$$\frac{1}{m}I_{m, ma}(\mF)\subset \frac{a}{ma}\kb_{ma}(\mF).$$
This implies $\lct(X, \Delta; a\cdot \kb_\bullet(\mF))\geq 1$ and hence $\lct(X, \Delta; \kb_\bullet(\mF))\geq a$.
\end{proof}

\begin{lemma}\label{lem: slope-delta}
Let $\mF$ be a filtration on $R$. Suppose $\mF^0R_m=R_m$ for any $m\in M(L)$. Then $\mu_{X, \Delta}(\mF)\geq \min\{1, \delta(X, \Delta)\}\cdot S_{X, \Delta}(\mF)$.
\end{lemma}

\begin{proof}
If $\delta(X, \Delta)\geq 1$, we are done by \cite[Corollary 4.5]{XZ20b}. So we just assume $\delta:=\delta(X, \Delta)<1$. Let $D\in \frac{1}{N}|-N(K_X+\Delta)|$ be a general element, where $N$ is a sufficiently divisible positive integer. By \cite[Theorem 1.8]{LXZ22}, $(X, \Delta+(1-\delta)D)$ is a K-semistable log Fano $\bQ$-pair (note that $\delta$ is rational by \cite[Theorem 1.2]{LXZ22}). The filtration $\mF$ naturally induces a filtration $\mF'$ on $R':=\oplus_{m\in M(L')}H^0(X, mL')$, where $L':=-K_X-\Delta-(1-\delta)D$. By \cite[Corollary 4.5]{XZ20b}, we have
$$\mu_{X, \Delta+(1-\delta)D}(\mF') \geq S_{X, \Delta+(1-\delta)D}(\mF').$$
An easy computation implies $S_{X, \Delta+(1-\delta)D}(\mF')=\delta\cdot S_{X, \Delta}(\mF)$. Thus it suffices to show that $\mu_{X, \Delta}(\mF)\geq \mu_{X, \Delta+(1-\delta)D}(\mF')$. For sufficiently divisible $m\in \bN$ and $\lambda>0$, we write
$$I_{m, \lambda}:= I_{m, \lambda}(\mF) \quad \text{and}\quad J_{m, \lambda}:= I_{m, \lambda}(\mF').$$
Then we see $J_{m, \lambda}=I_{m\delta, \lambda}$ and $\frac{1}{m}J_{m, mt}=\frac{1}{m} I_{m\delta, mt}\subset \frac{1}{m}I_{m, mt}$, which implies $J_\bullet^{(t)}\subset I_\bullet^{(t)}$. For any rational number $a< \mu_{X, \Delta+(1-\delta)D}(\mF')$, we have $\lct(X, \Delta+(1-\delta)D; J_\bullet^{(a)})\geq 1$. This implies $\lct(X, \Delta+(1-\delta)D; I_\bullet^{(a)})\geq 1$ and thus $\lct(X, \Delta; I_\bullet^{(a)})\geq 1$. Hence $\mu_{X, \Delta}(\mF)\geq a$. The proof is complete.
\end{proof}

\begin{lemma}\label{lem: ding-beta}
Let $\mF$ be a filtration on $R$. Then $\fD^{\NA}(\mF)\leq \beta_{X, \Delta}(\mF)$.
\end{lemma}
 \begin{proof}
 Proved by \cite[Theorem 4.4]{XZ20b}.
 \end{proof}

\subsection{Separatedness}

Let $\alpha: (X, \Delta)\to C$ and $\alpha': (X', \Delta')\to C$ be two $\bR$-Gorenstein families of log Fano pairs over a smooth pointed curve $0\in C$ such that
\begin{enumerate}
\item the two families are isomorphic over $C^\circ:=C\setminus 0$;
\item the central fibers $(X_{0}, \Delta_{0})$ and $(X'_{0}, \Delta'_{0})$ are both K-semistable.
\end{enumerate}
Up to shrinking $C$ around $0$, we may assume $C$ is affine and $\alpha, \alpha'$ are two families of K-semistable log Fano pairs (e.g. Theorem \ref{thm: openness}). In the following, we write 
$$\Delta=\sum_{j=1}^k a_j D_{j} \quad \text{and}\quad \Delta'=\sum_{j=1}^k a_j D'_{j}, $$ where $D_j$'s and $D'_j$'s are effective Weil divisors. For any vector $\fr:=(r_1,...,r_k)$, we write $\Delta(\fr):=\sum_{j=1}^k r_j D_j$ and $\Delta'(\fr):=\sum_{j=1}^k r_j D'_j$. Respectively, we write $\Delta(\fr)_0:=\sum_{j=1}^k r_j D_{j,0}$ and $\Delta'(\fr)_0:=\sum_{j=1}^k r_j D'_{j, 0}$ for the central fibers.

The following result on separatedness is well known by \cite{BX19} when $\fa:=(a_1,...,a_k)$ is rational, so we just assume $\fa$ is irrational.

\begin{theorem}\label{thm: R-BX}
The two central fibers $(X_{0}, \Delta_{0})$ and $(X'_{0}, \Delta'_{0})$ are $S$-equivalent to each other. In other words, they share a common K-semistable degeneration by two special test configurations.
\end{theorem}

\begin{proof}
We divide the proof into several steps.

\

\textbf{Step 1}. \textit{In this step, we construct two filtrations for a rational vector.} By Proposition \ref{prop: general QR}, there exists a rational polytope $P$ containing $\fa$ as an interior point such that $(X, \Delta(\fr))\to C$ and $(X', \Delta'(\fr))\to C$ are still $\bR$-Gorenstein families of log Fano pairs for any $\fr:=(r_1,...,r_k)\in P$. We may assume $P$ is of minimal dimension. \textit{Choose $\fr\in P$ to be rational.}
Let $r$ be a sufficiently divisible positive integer such that $rL(\fr):=-r(K_X+\Delta(\fr))$ and $rL'(\fr):=-r(K_{X'}+\Delta'(\fr))$ are very ample. Define
$$\mR_m(\fr):=H^0(X, mrL(\fr)),\ \ \ \mR_m'(\fr):=H^0(X', mrL'(\fr)), $$
$$R_m(\fr):=H^0(X_0, mrL(\fr)_0),\ \ \ R_m'(\fr):=H^0(X_0', mrL'(\fr)_0), $$
$$\mR(\fr):=\bigoplus_{m\in \bN}\mR_m(\fr),\ \ \ \mR'(\fr):=\bigoplus_{m\in \bN}\mR'_m(\fr),$$
$$R(\fr):=\bigoplus_{m\in \bN}R_m(\fr),\ \ \ R'(\fr):=\bigoplus_{m\in \bN}R'_m(\fr). $$
By Kawamata-Viehweg vanishing, we see that $R^i\alpha_*\mO_X(mrL(\fr))$ and $R^i\alpha'_*\mO_{X'}(mrL'(\fr))$ vanish for $i>0$ and $m> 0$. 
Hence the restriction maps
$$\mR_m(\fr)\otimes \bC(0)\to R_m(\fr) \quad \text{and} \quad \mR'_m(\fr)\otimes\bC(0) \to R_m'(\fr)$$
are both  isomorphic for $m\in \bN$,  and $\alpha_*\mO_X(mrL(\fr))$ and $\alpha'_*\mO_{X'}(mrL'(\fr))$ are both vector bundles for $m\in \bN$ and their cohomology commutes with base change.

Now we begin to construct two filtrations of $\mR(\fr)$ and $\mR'(\fr)$. We may assume that the central fibers $X_0$ and $X_0'$ are not birational to each other, otherwise $X$ and $X'$ are isomorphic in codimension one and $\mR(\fr)\cong \mR'(\fr)$, which implies that they are in fact isomorphic.
Note that $\ord_{X_0}$ (resp. $\ord_{X_0'}$) is a divisorial valuation on $X'$ (resp. $X$). We define the following filtrations for $m\in \bN$ and $p\in \bZ$:
$$\mF^p\mR_m(\fr):=\{s\in \mR_m(\fr)\ |\ \ord_{X_0'}(s)\geq p\},\ \ \  \mF'^p\mR'_m(\fr):=\{s\in \mR'_m(\fr)\ |\ \ord_{X_0}(s)\geq p\},$$
$$\mF^pR_m(\fr):={\rm Im}\{\mF^p\mR_m(\fr)\otimes \bC(0)\to R_m(\fr)\},$$ 
$$\mF'^pR'_m(\fr):={\rm Im}\{\mF'^p\mR'_m(\fr)\otimes \bC(0)\to R'_m(\fr)\} .$$
Up to shrinking $C$ around $0$,  there exists $\pi\in \mO(C)$ such that ${\rm{div}}_C(\pi)=[0]$.  It is clear that we have the following relations for $m\in  \bN$:
$$R_m(\fr)\cong \mR_m(\fr)/\pi\mR_m(\fr) \quad \text{and} \quad \mF^pR_m(\fr)\cong\frac{\mF^p\mR_m(\fr)}{\mF^p\mR_m(\fr)\cap \pi\mR_m(\fr)},$$
and the similar statement holds for $\mF'$. Therefore,  we see
$$\gr_\mF^pR_m(\fr):=\mF^pR_m(\fr)/\mF^{p+1}R_m(\fr)\cong \frac{\mF^p\mR_m(\fr)}{\mF^{p+1}\mR_m(\fr)+(\mF^p\mR_m(\fr)\cap \pi\mR_m(\fr))}. $$
Similarly, 
$$\gr_{\mF'}^pR'_m(\fr):=\mF'^pR'_m(\fr)/\mF'^{p+1}R'_m(\fr)\cong \frac{\mF'^p\mR'_m(\fr)}{\mF'^{p+1}\mR'_m(\fr)+(\mF'^p\mR'_m(\fr)\cap \pi\mR'_m(\fr))}. $$

\

\textbf{Step 2}. \textit{In this step, we study the relationship between the two filtrations $\mF$ and $\mF'$.} Let $\phi: Y\to X$ and $\phi': Y\to X'$ be the common resolution
\begin{center}
	\begin{tikzcd}[column sep = 2em, row sep = 2em]
	&Y\arrow[ld,"\phi",swap]\arrow[rd,"\phi'"]&\\
	 X && X' 
        \end{tikzcd}
\end{center}
and denote $V$ and $W$ for the birational transform of $X_0$ and $X_0'$ on $Y$ respectively. Then there are some rational numbers $a(\fr)$ and $a'(\fr)$ such that
$$\phi^*L(\fr)-\phi'^*L'(\fr)= a(\fr) W-a'(\fr) V+P,$$
where $P$ is both $\phi$-exceptional and $\phi'$-exceptional, and
$$a(\fr)=A_{X, \Delta(\fr)+X_0}(W) \quad \text{and} \quad a'(\fr)=A_{X', \Delta'(\fr)+X_0'}(V).$$
Therefore, 
\begin{align*}
\mF^p\mR_m(\fr)&\cong\ H^0(Y, \phi^*mrL(\fr)-pW)\\
&=\ H^0(Y, \phi'^*mrL'(\fr)-mra'(\fr)V+(mra(\fr)-p)W+mrP).
\end{align*}
For $s\in \mF^p\mR_m(\fr)$, multiplying $\phi^*s$ by $\pi^{mra(\fr)-p}$ gives an element in
$$\mF'^{mr(a(\fr)+a'(\fr))-p}\mR'_m(\fr)\cong H^0(Y, \phi'^*mrL'(\fr)-(mra(\fr)+mra'(\fr)-p)V). $$
Thus for each $p\in \bZ$ and $m\in \bN$, we have a map
$$\phi_{m,p}: \mF^p\mR_m(\fr)\to \mF'^{mr(a(\fr)+a'(\fr))-p}\mR'_m(\fr),\ \ s\mapsto \pi^{mra(\fr)-p}\cdot \phi^*s . $$
Similarly, we also have
$$\phi'_{m,p}: \mF'^p\mR'_m(\fr)\to \mF^{mr(a(\fr)+a'(\fr))-p}\mR_m(\fr),\ \ s'\mapsto \pi^{mra'(\fr)-p}\cdot \phi'^*(s'). $$
It is not hard to see by the construction that $\phi'_{m,mr(a(\fr)+a'(\fr))-p}\circ \phi_{m,p}=\id$, thus $\phi_{m,p}$ and $\phi'_{m,p}$ are isomorphisms for $p\in \bZ$ and $m\in \bN$.

\

\textbf{Step 3}. \textit{In this step, we show that $\mF$ and $\mF'$ induce two isomorphic graded filtrations on $R_m(\fr)$ and $R_m'(\fr)$.}
Recalling in step 1,  we have
$$\gr_\mF^pR_m(\fr):=\mF^pR_m(\fr)/\mF^{p+1}R_m(\fr)\cong \frac{\mF^p\mR_m(\fr)}{\mF^{p+1}\mR_m(\fr)+(\mF^p\mR_m(\fr)\cap \pi\mR_m(\fr))}, $$
and
$$\gr_{\mF'}^pR'_m(\fr):=\mF'^pR'_m(\fr)/\mF'^{p+1}R'_m(\fr)\cong \frac{\mF'^p\mR'_m(\fr)}{\mF'^{p+1}\mR'_m(\fr)+(\mF'^p\mR'_m(\fr)\cap \pi\mR'_m(\fr))}. $$
We first show that 
\begin{align*}
&\phi_{m,p}\left(\mF^{p+1}\mR_m(\fr)+(\mF^p\mR_m(\fr)\cap \pi\mR_m(\fr))\right)\\
\subset\  \ & \mF'^{mr(a(\fr)+a'(\fr))-p+1}\mR'_m(\fr)+(\mF^{mr(a(\fr)+a'(\fr))-p}\mR'_m(\fr)\cap \pi\mR_m'(\fr)).
\end{align*}
To see this, we first choose $s\in \mF^p\mR_m(\fr)\cap \pi\mR_m(\fr)$, then
\begin{align*}
& s\in H^0(Y, \phi^*mrL(\fr)-pW-V)\\
=\ \ &H^0(Y, \phi'^*mrL'(\fr)-(mra'(\fr)+1)V+(mra(\fr)-p)W+mrP). 
\end{align*}
Thus $\phi_{m,p}(s)\in \mF'^{mr(a(\fr)+a'(\fr))-p+1}\mR'_m(\fr)$. Next we choose $s\in \mF^{p+1}\mR_m(\fr)$. Similarly, we have $\phi_{m,p}(s)\in \mF'^{mr(a(\fr)+a'(\fr))-p}\mR'_m(\fr)$. On the other hand,  there is a rational function $f\in K(X)$ such that 
\begin{align*}
&\ord_{W}({\rm div}(f)+\phi^*mrL(\fr))\\
=\ \  &\ord_{W}({\rm div}(f)+\phi'^*mrL'(\fr)+mra(\fr)W-mra'(\fr)V+mrP) \\
\geq \ \  & p+1 .
\end{align*}
It implies the following inequality:
$$\ord_{X'_0}({\rm div}(f)+mrL'(\fr)) \geq p+1-mra(\fr). $$
Thus $\ord_{X_0'}(\pi^{mra(\fr)-p}\cdot f)\geq 1$ and $\phi_{m,p}(s)\in \pi\mR'_m(\fr)$, and hence
$$\phi_{m,p}(s)\in \mF'^{mr(a(\fr)+a'(\fr))-p}\mR'_m(\fr)\cap \pi\mR'_m(\fr) $$
In summary we obtain
\begin{align*}
&\phi_{m,p}\left(\mF^{p+1}\mR_m(\fr)+(\mF^p\mR_m(\fr)\cap \pi\mR_m(\fr))\right)\\
\subset\  \ & \mF'^{mr(a(\fr)+a'(\fr))-p+1}\mR'_m(\fr)+(\mF^{mr(a(\fr)+a'(\fr))-p}\mR'_m(\fr)\cap \pi\mR_m'(\fr)).
\end{align*}
Therefore, there is a morphism induced by $\{\phi_{m,p}\}_{m\in \bN, p\in \bZ}$:
$$\tilde{\phi}: \bigoplus_{m\in \bN}\bigoplus_{p\in \bZ}\gr^p_{\mF} R_m(\fr)\to \bigoplus_{m\in \bN}\bigoplus_{p\in \bZ}\gr^p_{\mF'} R'_m(\fr), $$
which sends $\gr^p_{\mF} R_m(\fr)$ to $\gr_{\mF'}^{mr(a(\fr)+a'(\fr))-p}R'_m(\fr)$. Similarly, we could establish the inverse map
$$\tilde{\phi}^{-1}: \bigoplus_{m\in \bN}\bigoplus_{p\in \bZ}\gr^p_{\mF'} R'_m(\fr)\to \bigoplus_{m\in \bN}\bigoplus_{p\in \bZ}\gr^p_{\mF} R_m(\fr), $$
which sends $\gr^p_{\mF'} R'_m(\fr)$ to $\gr_{\mF}^{mr(a(\fr)+a'(\fr))-p}R_m(\fr)$. It is not hard to see that $\tilde{\phi}$ and $\tilde{\phi}^{-1}$ are both isomorphisms since $\phi_{m,p}$ and $\phi'_{m,p}$ are both isomorphisms by step 2.

\

\textbf{Step 4}. \textit{In this step,  we study some invariants with respect  to the rational vector $\fr$.} To make a distinction, we denote $\mF\mR(\fr)$ (resp. $\mF R(\fr)$) to be the filtration on $\mR(\fr)$ (resp. $R(\fr)$). Similar notation applies to $\mF'$. Let  $\ka_{\bullet}(\mF\mR(\fr))$ be the base ideal sequence of $\mF\mR(\fr)$,  and let $\kb_{\bullet}(\mF R(\fr))$ be the restriction of $\ka_{\bullet}(\mF\mR(\fr))$ to $X_0$. Similarly, we have $\ka_{\bullet}(\mF'\mR'(\fr))$
and $\kb_{\bullet}(\mF' R'(\fr))$. By Lemma \ref{lem: slope-lct} and the inversion of adjunction we easily have the following $(\spadesuit)$:
\begin{align*}
&\mu(\mF R(\fr))\leq \lct\left(X_0, \Delta_0(\fr); \kb_\bullet(\mF R(\fr))\right)=\lct\left(X, \Delta(\fr)+X_0; \ka_\bullet(\mF\mR(\fr))\right)\leq a(\fr), \ \ \text{and} \\
&\mu(\mF' R'(\fr))\leq \lct\left(X'_0, \Delta'_0(\fr); \kb_\bullet(\mF' R'(\fr))\right)=\lct\left(X', \Delta'(\fr)+X'_0; \ka_\bullet(\mF'\mR'(\fr))\right) \leq a'(\fr).
\end{align*}
Denote $\delta(\fr):=\delta(X_0, \Delta(\fr)_0)$ and $\delta'(\fr):=\delta(X_0', \Delta'(\fr)_0)$. 
By Lemma \ref{lem: slope-delta}, we have $(\spadesuit\spadesuit)$
\begin{align*}
&\mu(\mF R(\fr)) \geq \min\{1, \delta(\fr)\}\cdot S(\mF R(\fr)),  \ \ \text{and} \\
&\mu(\mF' R'(\fr)) \geq \min\{1, \delta'(\fr)\}\cdot S(\mF' R'(\fr)).
\end{align*}
On the other hand, we have the following observation by step 3:
\begin{align*}
& \sum_{p=0}^{mr(a(\fr)+a'(\fr))}p\cdot \dim \gr_\mF^p R_m(\fr)+\sum_{p=0}^{mr(a(\fr)+a'(\fr))} p\cdot \dim \gr_{\mF'}^p R'_m(\fr)\\
=\ \ & \sum_{p=0}^{mr(a(\fr)+a'(\fr))}p\cdot \dim \gr_\mF^p R_m(\fr)+\sum_{p=0}^{mr(a(\fr)+a'(\fr))}p\cdot \dim \gr_{\mF}^{mr(a(\fr)+a'(\fr))-p}R_m(\fr)\\
=\ \ & mr(a(\fr)+a'(\fr))\sum_{p=0}^{mr(a(\fr)+a'(\fr))}\dim \gr_\mF^p R_m(\fr)\\
=\ \ & mr(a(\fr)+a'(\fr))\cdot \dim R_m(\fr).
\end{align*}
Let $m$ tend to $+\infty$, we have $(\clubsuit)$
$$a(\fr)+a'(\fr)=S(\mF R(\fr))+S(\mF' R'(\fr)). $$
By the continuity of $\delta$-invariants, we have
$\lim_{\fr\to \fa} \delta(\fr)\geq 1$ and $\lim_{\fr\to \fa}\delta'(\fr)\geq 1$ (note that $(X_0, \Delta_0)$ and $(X_0', \Delta_0')$ are K-semistable). Recalling the notation in step 2, we also have
$$\lim_{\fr\to \fa} a(\fr) =A_{X, \Delta+X_0}(W)\quad \text{and}\quad \lim_{\fr\to \fa} a'(\fr)=A_{X', \Delta'+X_0'}(V). $$
Combining $(\clubsuit)$ with $(\spadesuit)$ and $(\spadesuit\spadesuit)$, we see the following $(\clubsuit\clubsuit)$
\begin{align*}
&\lim_{\fr\to \fa} a(\fr)=\lim_{\fr\to \fa}\mu(\mF R(\fr))=\lim_{\fr\to \fa}\lct\left(X, \Delta(\fr)+X_0; \ka_\bullet(\mF \mR(\fr))\right), \ \ \text{and}\\
&\lim_{\fr\to \fa} a'(\fr)=\lim_{\fr\to \fa} \mu(\mF' R'(\fr))=\lim_{\fr\to \fa}\lct\left(X', \Delta'(\fr)+X'_0; \ka_\bullet(\mF' \mR'(\fr))\right).
\end{align*}

\

\textbf{Step 5}. \textit{In this step, we show that the following graded algebra is finitely generated for $\fr$ sufficiently close to $\fa$: }
$$\bigoplus_{m\in \bN}\bigoplus_{p\in \bZ}\gr^p_{\mF} R_m(\fr)\cong \bigoplus_{m\in \bN}\bigoplus_{p\in \bZ}\gr^p_{\mF'} R'_m(\fr). $$
By step 4, for any $0<\epsilon\ll 1$, we have $\mu(\mF R(\fr))>a(\fr)-\epsilon$ if $\fr\in P$ is sufficiently close to $\fa$.
This means that for $m\gg 1$, there exists a general element $D\in \mF^{mr(a(\fr)-\epsilon)} \mR_m(\fr)$ such that $(X, \Delta(\fr)+X_0+\frac{1}{mr}D)$ is log canonical with
$$A_{X, \Delta(\fr)+X_0+\frac{1}{mr}D}(X_0')\leq a(\fr)-\frac{1}{mr}\cdot \ord_{X_0'}(D)\leq  \epsilon.$$
By \cite{BCHM10}, there exists a birational model $\mu: Z\to X$ which precisely extracts $X_0'$, and the following graded algebra is finitely generated:
$$\bigoplus_{m\in \bN}\bigoplus_{p\in \bZ}\mF^p\mR_m(\fr)=\bigoplus_{m\in \bN}\bigoplus_{p\in \bZ} H^0(Z, \mu^*mrL(\fr)-p\widetilde{X_0'}),$$
where $\widetilde{X_0'}$ is the birational transform of $X_0'$ on $Z$.
By step 1, it is clear to see that $\bigoplus_{m\in \bN}\bigoplus_{p\in \bZ}\gr^p_{\mF} R_m(\fr)$ is also finitely generated as it is a quotient of the graded algebra $\bigoplus_{m\in \bN}\bigoplus_{p\in \bZ}\mF^p\mR_m(\fr)$.

\

\textbf{Step 6}. \textit{In this step, we construct a family over a two-dimensional base which encodes the common degeneration of $(X_0, \Delta(\fr)_0)$ and $(X_0', \Delta'(\fr)_0)$.} 
We may write $C=\Spec A$ and recall in step 1 that $\pi$ is a regular function on $C$ such that ${\rm{div}}(\pi)=[0]\in C$. Define
$$S:=\Spec A[s, t]/(st-\pi)\quad \text{and} \quad S^\circ:=S\setminus 0_{st},$$
where $0_{st}\in S$ is the closed point on $S$ defined by $s=t=0$.
Denote $0_s\subset S$ (resp. $0_t\subset S$) to be the divisor defined by $s=0$ (resp. $t=0$).
Applying the construction in \cite[Section 3]{ABHLX20}, we could glue $X\to C$ and $X'\to C$ along the isomorphism 
$$(X, L(\fr))\times_C C^\circ\cong (X', L'(\fr))\times_CC^\circ$$ 
to get a $\bG_m$-equivariant polarized family over $S^\circ$, denoted by
$$\Phi_\fr^\circ: (\mX^\circ, \mB(\fr)^\circ; \mL(\fr)^\circ)\to S^\circ, $$
where $\mB(\fr)^\circ:=\sum_{j=1}^kr_j \mD_j(\fr)^\circ$ is the extension of $\Delta(\fr)=\sum_{j=1}^kr_j D_j$, and $\mL(\fr)^\circ:=-K_{\mX^\circ/S^\circ}-\mB(\fr)^\circ$. By the finite generation obtained in step 5 and the construction in \cite[Section 3]{ABHLX20}, there exists an extension of $\Phi_\fr^\circ$, denoted by
$$\Phi_\fr: (\fX, \fB(\fr)+\fX_{0_s}+\fX_{0_t}; \fL(\fr))\to S, $$
where $\fB(\fr):=\sum_{j=1}^kr_j\fD_j$ is the natural extension of $\mB(\fr)^\circ$, $\fL(\fr)$ is a polarization on $\fX$ extending $\mL(\fr)^\circ$, and $\fX_{0_s}$ (resp. $\fX_{0_t}$) is the restriction of $\fX$ over $0_s$ (resp. $0_t$). From the construction (e.g. \cite[Section 3]{ABHLX20}), we have
\begin{enumerate}
\item the restriction of $(\fX, \fB(\fr)+\fX_{0_s}+\fX_{0_t}; \fL(\fr))\to S$ to $S\setminus (0_s\cup 0_t)$ is $\bG_m$-equivariantly isomorphic to $(X, \Delta(\fr), L(\fr))\times_C C^\circ\times \bG_m\to C^\circ\times \bG_m $;

\item $\fX_{0_s}\to \bA^1_t$ is isomorphic to 
$\Proj\bigoplus_{m\in \bN}\bigoplus_{p\in \bZ}\mF^pR_m(\fr)\cdot t^{-p}\to \bA_t^1;$

\item  $\fX_{0_t}\to \bA^1_s$ is isomorphic to 
$\Proj\bigoplus_{m\in \bN}\bigoplus_{p\in \bZ}\mF'^pR'_m(\fr)\cdot s^{-p}\to \bA_s^1;$

\item the central fiber $\fX_{0_{st}}$ (i.e. fiber over $0_{st}\in S$) of $\Phi_\fr$ is
$$\Proj\bigoplus_{m\in \bN}\bigoplus_{p\in \bZ}\gr^p_{\mF} R_m(\fr)\cong \Proj \bigoplus_{m\in \bN}\bigoplus_{p\in \bZ}\gr^p_{\mF'} R'_m(\fr).$$
\end{enumerate}

We emphasize here that $\fX$ constructed above depends on the choice of the rational vector $\fr\in P$, but we still use the notation $\fX$ for simplicity rather than $\fX_\fr$ as there should be no confusion. In the next step, we will remove this dependence up to shrinking $P$ around $\fa$.

\

\textbf{Step 7}. \textit{In this step, we show that the extending family $\Phi_\fr$ constructed in step 6 has good properties. More precisely, we show that the ambient space of $\Phi_\fr$ does not depend on $\fr$ up to shrinking $P$ around $\fa$. Moreover, for the family $\Phi_\fr: (\fX, \fB(\fr)+\fX_{0_s}+\fX_{0_t})\to S$, we aim to show that $\fL(\fr)=-K_{\fX/S}-\fB(\fr)$, and $(\fX, \fB(\fa)+\fX_{0_s}+\fX_{0_t})$ is log canonical (note that we replace $\fr$ with $\fa$).}

Similarly, as in step 6, by the construction (e.g. \cite[Section 3]{ABHLX20}), we can glue 
 $(X, \Delta(\fr); L)\to C$ and $(X, \Delta(\fr); L)\to C$ along the identity 
$$(X, \Delta(\fr); L(\fr))\times_C C^\circ\stackrel{\id}{\to} (X, \Delta(\fr); L(\fr))\times_C C^\circ$$ 
to obtain another family over $S$, denoted by $\hat{\Phi}_\fr: (\hat{\fX}, \hat{\fB}(\fr)+\hat{\fX}_{0_s}+\hat{\fX}_{0_t}; \hat{\fL}(\fr))\to S,$
such that the following conditions are satisfied:
\begin{enumerate}
\item the restriction $(\hat{\fX}_{0_s}, \hat{\fB}(\fr)_{0_s}; \hat{\fL}(\fr)_{0_s})\to \bA^1_t$ is isomorphic to $(X_0, \Delta(\fr)_0; L(\fr)_0)\times \bA_t^1\to \bA_t^1$;

\item the restriction $(\hat{\fX}_{0_t}, \hat{\fB}(\fr)_{0_t}; \hat{\fL}(\fr)_{0_t})\to \bA^1_s$ is isomorphic to $(X_0, \Delta(\fr)_0; L(\fr)_0)\times \bA_s^1\to \bA_s^1$;

\item $(\hat{\fX}, \hat{\fB}(\fr)+\hat{\fX}_{0_t}+\hat{\fX}_{0_s})$ is log canonical.
\end{enumerate}
Recall in step 5, $(X, \Delta(\fr)+X_0+\frac{1}{mr}D)$ is log canonical with
$$A_{X, \Delta(\fr)+X_0+\frac{1}{mr}D}(X_0')\leq a(\fr)-\frac{1}{mr}\cdot \ord_{X_0'}(D)\leq  \epsilon.$$
By \cite{BCHM10}, there exists a birational model $\mu: Z\to X$ which precisely extracts $X_0'$, and we write
$$K_Z+\widetilde{\Delta(\fr)}+ \widetilde{X_0}+\frac{1}{mr}\widetilde{D}+(1-\epsilon') \widetilde{X_0'}=\mu^*(K_X+\Delta(\fr)+X_0+\frac{1}{mr}D),$$
where $0\leq \epsilon'\leq \epsilon$. Let $\fD$ (resp. $\hat{\fD}$) be the extension of $D$ on $\fX$ (resp. $\hat{\fX}$), we clearly have
$$A_{\hat{\fX}, \hat{\fB}(\fr)+\hat{\fX}_{0_t}+\hat{\fX}_{0_s}+\frac{1}{mr}\hat{\fD}}(\fX_{0_t})=A_{X, \Delta(\fr)+X_0+\frac{1}{mr}D}(X_0')=1-\epsilon' .$$
By \cite{BCHM10}, there exists a crepant morphism which precisely extracts $\fX_{0_t}$:
$$\rho: (\widetilde{\fX}, \rho_*^{-1}(\hat{\fB}(\fr)+\hat{\fX}_{0_t}+\hat{\fX}_{0_s}+\frac{1}{mr}\hat{\fD})+(1-\epsilon')\widetilde{\fX_{0_t}})\to (\hat{\fX}, \hat{\fB}(\fr)+\hat{\fX}_{0_t}+\hat{\fX}_{0_s}+\frac{1}{mr}\hat{\fD}) ,$$
where $\widetilde{\fX_{0_t}}$ is the birational transform of $\fX_{0_t}$. Thus there is a birational contraction
$$(\widetilde{\fX}, \rho_*^{-1}(\hat{\fB}(\fr)+\hat{\fX}_{0_t}+\hat{\fX}_{0_s}+\frac{1}{mr}\hat{\fD})+(1-\epsilon')\widetilde{\fX_{0_t}})\dashrightarrow (\fX, \fB(\fr)+\fX_{0_s}+\frac{1}{mr}\fD+(1-\epsilon')\fX_{0_t}), $$
which precisely contracts $\rho_*^{-1}\hat{\fX}_{0_t}$. Since $(X, \Delta(\fr)+X_0+\frac{1}{mr}D)$ is log canonical Calabi-Yau (lc CY) over $C$, we see that $(\hat{\fX}, \hat{\fB}(\fr)+\hat{\fX}_{0_t}+\hat{\fX}_{0_s}+\frac{1}{mr}\hat{\fD})$ is lc CY over $S$.
This implies $(\fX, \fB(\fr)+\fX_{0_s}+\frac{1}{mr}\fD+(1-\epsilon')\fX_{0_t})$ is lc CY over $S$ and $\fX$ is Fano type over $S$ (in particular, $\fX$ is normal). 
Let 
$$\psi: (\fY, \psi_*^{-1} \fB(\fa))\to (\fX, \fB(\fa))$$ 
be a $\bG_m$-equivariant small $\bQ$-factorization and run $\bG_m$-equivariant MMP/$S$ on $-K_{\fY/S}-\psi_*^{-1} \fB(\fa)$ resulting in the anti-canonical model over $S$ with respect to $-K_{\fY/S}-\psi_*^{-1} \fB(\fa)$, denoted by 
$$\xi: (\fY, \psi_*^{-1} \fB(\fa))\dashrightarrow (\fZ, \xi_*\psi_*^{-1} \fB(\fa)).$$ 
By \cite[Corollary 1.1.5]{BCHM10}, there exists a rational polytope $P_1\subset P$ containing $\fa$ as an interior point such that $(\fZ, \xi_*\psi_*^{-1} \fB(\fr'))$ is the anti-canonical model of $(\fY, \psi_*^{-1} \fB(\fr'))$ for any rational vector $\fr'\in P_1$.  It is clear that $(\fZ, \xi_*\psi_*^{-1} \fB(\fr'))$ with its anti-canonical polarization is actually an extension of 
$$\Phi^\circ_{\fr'}: (\mX^\circ, \mB(\fr')^\circ; \mL(\fr')^\circ)\to S^\circ.$$
By the uniqueness of the extension (e.g. \cite[Lemma 2.16]{ABHLX20}), we see that the family $\Phi_\fr$ coincides with $(\fZ, \xi_*\psi_*^{-1} \fB(\fr))$ when $\fr\in P_1$, and we indeed have $\fL(\fr)=-K_{\fX/S}-\fB(\fr)$. This means that the ambient space for $\Phi_\fr$ actually does not depend on the choice of the rational vector $\fr\in P_1$.

On the other hand, we already see that $(\fX, \fB(\fr)+\fX_{0_s}+\frac{1}{mr}\fD+(1-\epsilon')\fX_{0_t})$ is log canonical. 
Since we could choose $\epsilon'$ arbitrarily small by choosing $\fr$ sufficiently close to $\fa$ (see step 5), we conclude that $(\fX, \fB(\fa)+\fX_{0_s}+\fX_{0_t})$ is also log canonical. 

\

\textbf{Step 8}. \textit{In this step, we make a summary and finish the proof}. In step 1, for each rational vector $\fr$ contained in a rational polytope $P$, we construct two filtrations for the graded rings associated to $(X, \Delta(\fr))$ and $(X', \Delta'(\fr))$.  In step 2-5, we show that the graded rings associated to the two filtrations are finitely generated (up to shrinking $P$), where we indeed rely on the K-semistability of $(X_0, \Delta_0)$ and $(X'_0, \Delta'_0)$. Based on the finite generation, in step 6, we construct an extending family $\Phi_\fr: (\fX, \fB(\fr)+\fX_{0_s}+\fX_{0_t}; \fL(\fr))\to S$ for each rational vector $\fr\in P$. In step 7, via MMP, we show that the extending family $\Phi_\fr$ actually does not depend on $\fr\in P$ (up to shrinking $P$ again), and $(\fX, \fB(\fa)+\fX_{0_s}+\fX_{0_t})$ is log canonical.

Note that $\fa$ is irrational. By the rationality of the log canonical polytope, we see that $(\fX, \fB(\fr)+\fX_{0_s}+\fX_{0_t})$ is log canonical for any rational vector $\fr\in P$ up to shrinking $P$. 
Applying adjunction, the central fiber of $(\fX, \fB(\fr)+\fX_{0_s}+\fX_{0_t})\to S$, i.e. $(\fX_{0_{st}}, \fB(\fr)_{0_{st}})$, admits slc singularities. In summary, we conclude that there exists a rational polytope $Q\subset P$ containing $\fa$ as an interior point such that for any $\fa'\in Q$, the following conditions are satisfied: 
\begin{enumerate}
\item the restriction of $(\fX, \fB(\fa'))\to S$ to $0_s$ gives a weakly special test configuration of $(X_0, \Delta(\fa')_0)$;
\item the restriction of $(\fX, \fB(\fa'))\to S$ to $0_t$ gives a weakly special test configuration of $(X'_0, \Delta'(\fa')_0)$;
\item the common degeneration of $(X_0, \Delta(\fa')_0)$ and $(X'_0, \Delta'(\fa')_0)$ is $(\fX_{0_{st}}, \fB(\fa')_{0_{st}})$.
\end{enumerate}
It remains to show that $(\fX_{0_s}, \fB(\fa)_{0_s})$ is a special test configuration of $(X_0, \Delta_0)$ such that the central fiber $(\fX_{0_{st}}, \fB(\fa)_{0_{st}})$ is also K-semistable. For any rational $\fa'\in Q$, by \cite[Theorem 3.2]{Fuj18}, we first have the following equality $(\star)$:
$$\Fut(\fX_{0_s}, \fB(\fa')_{0_s})=\Ding (\fX_{0_s}, \fB(\fa')_{0_s}).$$
By the construction in step 6, the weakly special test configuration $(\fX_{0_s}, \fB(\fa')_{0_s})$ is given by the filtration $\mF R(\fa')$. Applying the computation in \cite{Fuj18}, we have the following equality $(\star\star)$:
$$\Ding(\fX_{0_s}, \fB(\fa')_{0_s})=\fD^\NA(\mF R(\fa')). $$
Combining $(\star)$ and $(\star\star)$, by Lemma \ref{lem: ding-beta}, we have
$$\Fut(\fX_{0_s}, \fB(\fa')_{0_s})\leq \mu(\mF R(\fa'))-S(\mF R(\fa')). $$
Letting $\fa'$ tend to $\fa$, the right hand side of the above inequality tends to $0$ by step 4 (see $(\spadesuit)$, $(\spadesuit\spadesuit)$, $(\clubsuit)$). Thus we have $\Fut(\fX_{0_s}, \fB(\fa)_{0_s})\leq 0$. Since $(X_0, \Delta_0)$ is K-semistable, by Remark \ref{rem: weak LX}, we see $\Fut(\fX_{0_s}, \fB(\fa)_{0_s})= 0$. By Proposition \ref{prop: weak LX}, $(\fX_{0_s}, \fB(\fa)_{0_s})$ is indeed a special test configuration with a K-semistable central fiber.  The proof is complete.
\end{proof}

\subsection{$S$-completeness}\label{subsec: S}

In this subsection, we show that $\mM^\Kss_{d, v_0, \fa}$ (see Section \ref{subsec: stack}) is $S$-complete. Fix $R$ to be a DVR essentially of finite type over $\bC$ with fraction field $K$ and residue field $\kappa$.  Let $\pi\in R$ be a uniformizing parameter. Define the Artin stack 
$$\overline{\ST}_R:=[\Spec\left(R[s,t]/(st-\pi)\right)/\bC^*], $$
where $s$ and $t$ have weights $1$ and $-1$. Denote by $0\in \overline{\ST}_R$ the unique closed point defined by $s=t=0$. Note that $\overline{\ST}_R\setminus 0$ is the non-separated union $\Spec R \cup_{\Spec K} \Spec R$.

\begin{definition}
We say $\mM^\Kss_{d, v_0, \fa}$ is $S$-complete if any map $\overline{\ST}_R\setminus 0 \to \mM^\Kss_{d, v_0, \fa}$ extends uniquely to a map $\overline{\ST}_R \to \mM^\Kss_{d, v_0, \fa}$.
\end{definition}

Giving a map $\overline{\ST}_R\setminus 0 \to \mM^\Kss_{d, v_0, \fa}$ is equivalent to giving two families $[(X, \Delta)\to \Spec R]\in \mM^\Kss_{d, v_0, \fa}(R)$ and $[(X', \Delta')\to \Spec R]\in \mM^\Kss_{d, v_0, \fa}(R)$
such that they are isomorphic over $\Spec K$ (see the construction in \cite[Section 3]{ABHLX20}).

\begin{theorem}\label{thm: S}
The Artin stack $\mM^\Kss_{d, v_0, \fa}$ is $S$-complete.
\end{theorem}

\begin{proof}
When $\fa$ is rational, the result is well known by \cite{ABHLX20}. We just assume $\fa$ is irrational. Given two families $[(X, \Delta)\to \Spec R]\in \mM^\Kss_{d, v_0, \fa}(R)$ and $[(X', \Delta')\to \Spec R]\in \mM^\Kss_{d, v_0, \fa}(R)$ which are isomorphic over $\Spec K$.  Denote $$S:=\Spec(R[s,t]/(st-\pi)) \quad \text{and}\quad S^\circ:=S\setminus 0,$$
where $0\in S$ is the closed point defined by $s=t=0$ (note that we use $0_{st}$ in the proof of Theorem \ref{thm: R-BX}). The two families naturally glue to a family over $S^\circ$, denoted by
$$\Phi^\circ: (\mX^\circ, \mB^\circ; \mL^\circ)\to S^\circ, $$
where $\mL^\circ:=-K_{\mX^\circ/S}-\mB^\circ$ is a polarization.
By the proof of Theorem \ref{thm: R-BX}, there exists an extending family
$$\Phi: (\fX, \fB; \fL)\to S, $$
where $\fL:=-K_{\fX/S}-\fB$ is a polarization and the fiber over $0\in S$ is a K-semistable log Fano pair. By a similar argument as in the proof of Theorem \ref{thm: Theta} (see the step 4 of the proof), the extension is unique. In summary we see that any map $\overline{\ST}_R\setminus 0 \to \mM^\Kss_{d, v_0, \fa}$ indeed extends uniquely to a map $\overline{\ST}_R \to \mM^\Kss_{d, v_0, \fa}$. The proof is complete.
\end{proof}

\section{Properness of K-moduli}\label{sec: proper}

In this section, we aim to show that the Artin stack $\mM^{\Kss}_{d, v_0, \fa}$ (see Section \ref{subsec: stack}) admits a proper and separated good moduli space. We will avoid rebuilding the whole program of $\Theta$-stratifications as in $\bQ$-coefficients case (e.g. \cite{BHLLX21}), but try to reduce the problem to $\bQ$-coefficients case and conduct an approximation process.

\subsection{Properness}\label{subsec: proper}

Recall that we say a flat family of log Fano pairs $(X, \Delta)\to T$ over a smooth base $T$ is an $\bR$-Gorenstein family if $-K_{X/S}-\Delta$ is $\bR$-Cartier.  If every fiber is a log Fano $\bQ$-pair and $-K_{X/S}-\Delta$ is $\bQ$-Cartier, we say it is a \textit{$\bQ$-Gorenstein family}.
We denote $C^\circ=C\setminus 0$ to be a punctured smooth curve. 

\begin{proposition}\label{prop: proper Q}
Let $\pi: (X^\circ, \Delta^\circ)\to C^\circ$ be a $\bQ$-Gorenstein family of log Fano $\bQ$-pairs such that $\delta(X_t, \Delta_t)=\delta_0\leq 1$ for any $t\in C^\circ$. Then up to a finite base change, one could fill the missing fiber to get a complete $\bQ$-Gorenstein family of log Fano $\bQ$-pairs $(X, \Delta)\to C$ such that $\delta(X_0, \Delta_0)=\delta_0$.
\end{proposition}

When $\delta_0=1$, the result is well-known due to the properness of the Artin stack parametrizing K-semistable log Fano $\bQ$-pairs with fixed invariants (e.g. \cite{BHLLX21, LXZ22}). For the case $\delta_0<1$, the result is also well-known due to the existence of a well-ordered $\Theta$-stratification  on the stack of log Fano $\bQ$-pairs (see \cite[Definition 2.12]{BHLLX21} or \cite[Definition 6.1]{AHLH23}).
However, we prefer to give a more direct proof here relying on the case for $\delta_0=1$.

\begin{proof}[Proof of Proposition \ref{prop: proper Q}]
We only assume $\delta_0<1$. 
Fix a closed point $t\in C^\circ$ and consider the linear system $\frac{1}{m}|-m(K_{X^\circ_t}+\Delta^\circ_t)|$ for a sufficiently divisible $m\in \bN$. By \cite[Theorem 5.4]{LXZ22}, there exists an element $D\in \frac{1}{m}|-m(K_{X^\circ_t}+\Delta^\circ_t)|$ such that $(X^\circ_t, \Delta^\circ_t+(1-\delta_0)D)$ is a K-semistable log Fano $\bQ$-pair with $\delta(X^\circ_t, \Delta^\circ_t+(1-\delta_0)D)=1$. We may assume $C^\circ$ is affine and consider the following exact sequence:
$$0\to \mO_{X^\circ}(-m(K_{X^\circ}+\Delta^\circ)-X_t^\circ)\to \mO_{X^\circ}(-m(K_{X^\circ}+\Delta^\circ))\to \mO_{X^\circ_t}(-m(K_{X_t^\circ}+\Delta_t^\circ))\to 0. $$
Applying $\pi_*$ to the above exact sequence, we see  
$$H^0(X^\circ, -m(K_{X^\circ}+\Delta^\circ))\to H^0(X^{\circ}_t, -m(K_{X^{\circ}_t}+\Delta^{\circ}_t))$$
is surjective since $H^1(X^\circ, -m(K_{X^\circ}+\Delta^\circ)-X_t^{\circ})=0$ for sufficiently large $m$ by Serre vanishing. This means there exists a divisor $\mD^\circ\sim_\bQ -K_{X^\circ}-\Delta^\circ$ on $X^\circ$ such that $\mD^\circ_t=D$. By \cite{BLX22, Xu20}, up to shrinking the base, we may assume $(X^\circ, \Delta^\circ+(1-\delta_0)\mD^\circ)\to C^\circ$ is a family of K-semistable log Fano pairs (with $\delta$-invariants being one). By the properness of K-moduli for the case $\delta_0=1$, up to a finite base change, there exists a filling $(X_0, \Delta_0+(1-\delta_0)\mD_0)$ which is K-semistable with $\mD_0\sim_\bQ -K_{X_0}-\Delta_0$. It is clear that we have $\delta(X_0, \Delta_0)\geq \delta_0$. By the lower semi-continuity of $\delta$-invariants (e.g. \cite{BL22}) we see that $\delta(X_0, \Delta_0)=\delta_0$. The proof is complete.
\end{proof}

We are ready to confirm the properness of $\mM^\Kss_{d, v_0, \fa}$.

\begin{theorem}\label{thm: proper R}
Let $\pi: (X^\circ, \Delta^\circ)\to C^\circ$ be an $\bR$-Gorenstein family of K-semistable log Fano $\bR$-pairs. Then up to a finite base change, one could fill the missing fiber to get a complete $\bR$-Gorenstein family of log Fano $\bR$-pairs $(X, \Delta)\to C$ such that $(X_0, \Delta_0)$ is K-semistable. In particular, $\mM^\Kss_{d, v_0, \fa}$ satisfies the existence part of valuative criterion for the properness.
\end{theorem}

\begin{proof}
We divide the proof into two parts: (1) there exists a fiber of $\pi$ with $\delta$-invariant $>1$; (2) every fiber of $\pi$ admits $\delta$-invariant $=1$.

\

Firstly, we assume there exists a closed point $t_0\in C^\circ$ such that $\delta(X^{\circ}_{t_0}, \Delta^{\circ}_{t_0})>1$. By Theorem \ref{thm: openness}, up to shrinking $C^\circ$, we may assume $\delta(X^\circ_t, \Delta^\circ_t)> 1$ for any $t\in C^\circ$.
Write $\Delta^\circ:=\sum_{j=1}^s r_j D_j^\circ$, where $D_j^\circ$'s are effective Weil divisors and $(r_1,...,r_s)$ is irrational. By Proposition \ref{prop: general QR}, there exists a rational polytope $P$ containing $(r_1,...,r_s)$ as an interior point such that $(X^\circ, \sum_{j=1}^sx_j D_j^\circ)\to C^\circ$ is a $\bQ$-Gorenstein family of log Fano $\bQ$-pairs for any rational vector $(x_1,...,x_s)\in P$. We may assume the dimension of $P$ is the same as the rank of $(r_1,...,r_s)$ over $\bQ$. It is clear that there exists a positive real number $\epsilon_0>0$ such that $\vol(-K_{X^\circ_t}-\sum_{j=1}^sx_jD^\circ_{j,t})\geq \epsilon_0$ for any $t\in C^\circ$ and any $(x_1,...,x_s)\in P$.
Let $\{(r_{i1},...,r_{is})\}_{i=1}^\infty\subset P$ be a sequence of rational vectors tending to $(r_1,...,r_s)$. By Theorem \ref{thm: openness} and the continuity of $\delta$-invariants,  for each $i\gg 1$, we may assume $(X^\circ, \sum_{j=1}^sr_{ij}D^\circ_j)\to C^\circ$ is a $\bQ$-Gorenstein family of K-semistable log Fano $\bQ$-pairs up to shrinking $C^\circ$. By the properness of the K-moduli of K-semistable log Fano $\bQ$-pairs (e.g. \cite{BHLLX21, LXZ22}), up to a finite base change, one could fill the punctured family $(X^\circ, \sum_{j=1}^sr_{ij}D_j^\circ)\to C^\circ$ to get a complete $\bQ$-Gorenstein family of K-semistable log Fano $\bQ$-pairs $(X^{(i)}, \sum_{j=1}^sr_{ij}D^{(i)}_{j})\to C$. Consider the set of fillings:
$$\mP_1:=\{(X^{(i)}_0, \sum_{j=1}^s D^{(i)}_{j,0} )\}_i.$$
Note that for each index $i$, we have $(X^{(i)}_0, \sum_{j=1}^s r_{ij}D^{(i)}_{j,0} )$ is K-semistable and $\vol(-K_{X^{(i)}_0}-\sum_{j=1}^s x_jD^{(i)}_{j,0})\geq \epsilon_0$ for any $(x_1,...,x_s)\in P$. By \cite[Proposition 7.1]{LZ24}, the set $\mP_1$ is log bounded. By \cite[Propositions 3.4, 3.5, 4.1]{LZ24}, up to choosing a countable subset of $\mP_1$, we may assume there exists a rational polytope $P_1\subset P$ containing $(r_1,...,r_s)$ as an interior point such that $\{(r_{i1},...,r_{is})\}_{i=1}^\infty\subset P_1$ and $(Y, \sum_{j=1}^s x_jB_j)$ is log Fano for any $(Y, \sum_{j=1}^s B_j)\in \mP_1$ and any $(x_1,...,x_s)\in P_1$. By \cite[Proposition 7.4]{LZ24}, up to shrinking $P_1$, 
we may assume
$$(X^{(i)}, \sum_{j=1}^sr_{j}D^{(i)}_{j})\to C$$
is an $\bR$-Gorenstein family of log Fano $\bR$-pairs for any $i$. It suffices to show that $(X^{(i)}_0, \sum_{j=1}^s r_jD^{(i)}_{j,0})$ is K-semistable for some $i$. 
Let $f: (\mX, \sum_{j=1}^s  \mD_j)\to T$ be a family of projective couples over a finite type scheme $T$, where $\mD_j$'s are effective Weil divisors, such that $\mP_1$ is contained in the set of fibers of $f$. By the proof of \cite[Theorem 5.4]{LZ24}, up to a base change $T'\to T$, we may assume the following conditions:
\begin{enumerate}
\item a countable subset of $\mP_1$ is contained in the set of fibers of $f': (\mX', \sum_{j=1}^s\mD_j')\to T'$ (obtained by the base change);
\item $(\mX'_t, \sum_{j=1}^s x_j\mD'_{j,t})$ is log Fano for any $t\in T'$ and any $(x_1,...,x_s)\in P_1$;
\item $\min\{\delta(\mX'_t, \sum_{j=1}^s x_j\mD'_{j,t}), 1\}$ does not depend on $t\in T'$ for any $(x_1,...,x_s)\in P_1$. 
\end{enumerate}
Up to subtracting redundant $i$, we may assume $\mP_1$ is contained in the set of fibers of $f'$. Thus for any $(X^{(i_0)}_0, \sum_{j=1}^s D^{(i_0)}_{j,0})\in \mP_1$, we see that $(X^{(i_0)}_0, \sum_{j=1}^s r_{ij}D^{(i_0)}_{j,0})$ is a K-semistable log Fano pair for any index $i$. By the continuity, we derive that $(X^{(i_0)}_0, \sum_{j=1}^s r_jD^{(i_0)}_{j,0})$ is K-semistable.
The proof for the first part is complete.

\

Now we turn to the second part and assume $\delta(X^\circ_t, \Delta^\circ_t)=1$ for every $t\in C^\circ$. The idea is the same as the first part with minor modification. We still use the notation in the first part. The different point is that for each $i$, $(X^\circ, \sum_{j=1}^sr_{ij}D^\circ_j)\to C^\circ$ is still a $\bQ$-Gorenstein family of  log Fano $\bQ$-pairs but not necessarily a family of K-semistable log Fano $\bQ$-pairs up to shrinking $C^\circ$. This, however, does not matter. By Theorem \ref{thm: constructible} and Theorem \ref{thm: openness}, up to shrinking $C^\circ$, we may assume every fiber of $(X^\circ, \sum_{j=1}^sr_{ij}D^\circ_j)\to C^\circ$ admits the same $\delta$-invariant, denoted by $\delta_i$. If $\delta_i\geq 1$, we could find the missing fiber as in the first part; if $\delta_i<1$, we could find the missing fiber by Proposition \ref{prop: proper Q}. In summary, for each index $i$, up to a finite base change and shrinking $C^\circ$, one could fill the punctured family $(X^\circ, \sum_{j=1}^sr_{ij}D_j^\circ)\to C^\circ$ to get a complete $\bQ$-Gorenstein family of log Fano $\bQ$-pairs $(X^{(i)}, \sum_{j=1}^sr_{ij}D^{(i)}_{j})\to C$ such that every fiber admits the same $\delta$-invariant $\delta_i$. It is clear to see $\delta_i$ tends to one.
Up to subtracting redundant $i$, we may assume $\delta_i$ is bounded from below by a positive real number $\delta_0$. As before, we consider the set of fillings:
$$\mP_1:=\{(X^{(i)}_0, \sum_{j=1}^s D^{(i)}_{j,0} )\}_i.$$
Thus for each $i$, we have $\delta(X^{(i)}_0, \sum_{j=1}^s r_{ij}D^{(i)}_{j,0} )\geq \delta_0$ and $\vol(-K_{X^{(i)}_0}-\sum_{j=1}^s x_jD^{(i)}_{j,0})\geq \epsilon_0$ for any $(x_1,...,x_s)\in P$. The rest of the proof is the same as part one, except that $(X^{(i_0)}_0, \sum_{j=1}^s r_{ij}D^{(i_0)}_{j,0})$ admitting $\delta$-invariant $\delta_i$ rather than being K-semistable.  Recall that $\delta_i$ tends to one. By the continuity, we still derive that $(X^{(i_0)}_0, \sum_{j=1}^s r_jD^{(i_0)}_{j,0})$ is K-semistable. The proof for the second part is complete.
\end{proof}

\subsection{Proper K-moduli space}\label{subsec: proper moduli}

Based on the preparations in Sections \ref{sec: stack}, \ref{sec: Theta}, \ref{sec: S}, \ref{subsec: proper}, 
we are ready to show in this subsection that the Artin stack $\mM^{\Kss}_{d, v_0, \fa}$ (see Section \ref{subsec: stack}) indeed admits a proper good moduli space.  

\begin{definition}\label{def: gms}
An algebraic space $Y$ is called a good moduli space of an Artin
stack $\mY$, if there is a quasi-compact morphism $\phi: \mY\to Y$ such that
\begin{enumerate}
\item  $\phi_*$ is an exact functor on quasi-coherent sheaves; and
\item $\phi_*\mO_\mY=\mO_Y$
\end{enumerate}
\end{definition}

\begin{theorem}{\rm (\cite[Theorem A]{AHLH23})}\label{thm: AHLH}
Let $\mY$ be an Artin stack of finite type with affine diagonal over
$\bC$, then $\mY$ admits a separated good moduli space if $\mY$ is $S$-complete and $\Theta$-reductive. If $\mY$ in addition satisfies the valuative criterion for properness, then $\mY$ admits a proper moduli space.
\end{theorem}

We are ready to construct the good moduli space of $\mM^\Kss_{d, v_0, \fa}$ (see Section \ref{subsec: stack}).

\begin{theorem}\label{thm: gms}
The Artin stack $\mM^\Kss_{d, v_0, \fa}$ admits a separated and proper good moduli space parametrizing K-polystable objects.
\end{theorem}

\begin{proof}
The result is well known when $\fa$ is rational (e.g. \cite{LXZ22}). We just assume $\fa$ is irrational. By the proof of Theorem \ref{thm: artin stack}, $\mM^\Kss_{d, v_0, \fa}$ automatically admits affine diagonal. 
By Theorem \ref{thm: Theta} and Theorem \ref{thm: S}, $\mM^\Kss_{d, v_0, \fa}$ is $\Theta$-reductive and $S$-complete. By Theorem \ref{thm: proper R}, $\mM^\Kss_{d, v_0, \fa}$ satisfies the valuative criterion for properness. Applying Theorem \ref{thm: AHLH}, we see that $\mM^\Kss_{d, v_0, \fa}$ admits a separated and proper good moduli space parametrizing $S$-equivalent classes. By Theorem \ref{thm: kps deg}, each $S$-equivalent class corresponds to a K-polystable log Fano $\bR$-pair. The proof is complete.
\end{proof}

\begin{corollary}\label{cor: reductivity}
Let $(X, \Delta)$ be a K-polystable log Fano $\bR$-pair. Then $\Aut(X, \Delta)$ is reductive.
\end{corollary}
\begin{proof}
See \cite[Proposition 12.14]{Alper13} or \cite[Theorem 8.5]{Xu24}.
\end{proof}

\bibliography{reference.bib}
\end{document}